\documentclass[11pt]{article}
\usepackage[margin=1in]{geometry}
\usepackage{amsmath,amssymb,amsthm,mathtools}
\usepackage{enumitem,algorithm,float,algpseudocode,microtype,xcolor,etoolbox,booktabs,graphicx}
\usepackage[numbers,square,sort&compress]{natbib}
\usepackage[colorlinks=true,linkcolor=blue,citecolor=blue,urlcolor=blue]{hyperref}
\usepackage{authblk}
\allowdisplaybreaks[1]
\newcommand{\R}{\mathbb R}
\newcommand{\E}{\mathbb E}
\newcommand{\tr}{\operatorname{tr}}
\newcommand{\rank}{\operatorname{rank}}
\newcommand{\range}{\operatorname{range}}
\newcommand{\vecop}{\operatorname{vec}}
\newcommand{\matop}{\operatorname{mat}}
\newcommand{\argmin}{\operatorname*{arg\,min}}
\newcommand{\argmax}{\operatorname*{arg\,max}}

\newcommand{\diag}{\operatorname{diag}}

\newtheorem{theorem}{Theorem}[section]
\newtheorem{proposition}[theorem]{Proposition}
\newtheorem{lemma}[theorem]{Lemma}
\newtheorem{corollary}[theorem]{Corollary}
\theoremstyle{definition}
\newtheorem{assumption}[theorem]{Assumption}
\newtheorem{example}[theorem]{Example}
\theoremstyle{remark}
\newtheorem{remark}[theorem]{Remark}

\numberwithin{equation}{section}

\newcommand{\PP}{\mathbb P}
\newcommand{\QQ}{\mathbb Q}

\ifdefined\subfigure\else
\newenvironment{subfigure}[1]{\begin{minipage}[t]{#1}\centering}{\end{minipage}}
\fi

\title{Distributionally Robust Risk Evaluation with \\Low-Rank Perturbations}
\author[1]{Linglingzhi Zhu\thanks{\href{mailto:llzzhu@gatech.edu}{llzzhu@gatech.edu}}}
\author[1]{Dongze Wu\thanks{\href{mailto:dwu381@gatech.edu}{dwu381@gatech.edu}}}
\author[1]{Yao Xie\thanks{\href{mailto:yao.xie@isye.gatech.edu}{yao.xie@isye.gatech.edu}}}
\affil[1]{H. Milton Stewart School of Industrial and Systems Engineering\authorcr Georgia Institute of Technology}

\date{November 2026}

\begin{document}
\maketitle
\begin{abstract}
We study distributionally robust risk under perturbations restricted to low-rank matrices whose row and column spaces are shared across observations, focusing on how much worst-case risk this restriction sacrifices and how to compute the restricted worst case. Under moment or Wasserstein budgets, we derive quantitative bounds on the gap between the unrestricted and low-rank worst-case risks, characterize when the gap vanishes, and identify the smallest lossless row and column dimensions. To generate worst-case perturbations, we apply projected gradient ascent to a curvature-corrected concave problem and derive computable risk certificates. We establish global sublinear convergence for the corrected objective and linear convergence to a worst-case distribution when the budget constraint is sufficiently active.
\end{abstract}



\section{Introduction}\label{sec:intro}
Distributionally robust optimization (DRO) chooses a decision by minimizing the largest expected loss over a prescribed class of probability distributions. Given a nominal distribution $\PP$ for an observation $\xi=(X,Y)$ and a loss function $\ell(\theta;\xi)$ that a decision $\theta$ incurs on the observation $\xi$, the problem takes the form
\[
    \inf_\theta\;\sup_{\QQ\in\mathcal A(\PP)} \E_{\xi'\sim \QQ}[\ell(\theta;\xi')],
\]
where $\mathcal A(\PP)$ is an ambiguity set describing the departures
from $\PP$ against which the decision is evaluated. 
Reweighting ambiguity sets, such as likelihood-ratio or $\varphi$-divergence balls with $\QQ\ll \PP$, are commonly used to modify the probabilities of nominal outcomes while leaving the outcomes themselves unchanged \citep{benTal2013,duchiNamkoong2021,gui2026distributionally}. Our focus instead is on uncertainty that acts directly on the input values, which can be modeled by ambiguity sets that allow the nominal outcomes themselves to move, as in moment-based and transportation-based formulations \citep{delage2010,xu2018,mohajerin2018wasserstein,blanchetMurthy2019,gaoKleywegt2023}. We therefore adopt a displacement view in which the input $X$ is perturbed to $X+Z$ and a budget controls the admissible perturbations $Z$.

For matrix-valued inputs $X\in\R^{m_1\times m_2}$, uncertainty can
concern not only the size of a change but also its row and column
structure. Many inputs in operations and engineering are matrices whose rows and columns
index different things: regions and time periods in power-system operations,
sensors and time in monitoring networks, or the two spatial axes of an image.
For such inputs a distributional shift has a direction as well as a size, and
realistic shifts can often be represented by a few patterns along each axis. We therefore
restrict the displacements to a row space and a column space shared by all observations,
\begin{equation}\label{eq:intro-model}
 Z=UCV^\top,\qquad C\mid\xi\sim q(\,\cdot\mid\xi),
\end{equation}
where $U\in\R^{m_1\times r}$ and $V\in\R^{m_2\times s}$ have orthonormal
columns, $C\in\R^{r\times s}$ is a random core matrix encoding the coordinates of the perturbation in the shared row and column bases, and $q$ is its conditional distribution.
We measure displacements by one of two quadratic budgets,
\begin{equation}\label{eq:intro-budgets}
 \E_q[\vecop(Z)\vecop(Z)^\top]\preceq B,\quad B\succ0,
 \qquad\text{or}\qquad
 \E_q\|Z\|_F^2\le\rho^2,\quad\rho>0,
\end{equation}
where $\E_q$ denotes expectation under $\PP(d\xi)\,q(dC\mid\xi)$. The matrix
second-moment budget bounds the mean-square size of every linear combination
of entries, with a matrix $B$ that can encode correlations across rows and
columns. The energy budget is a squared Wasserstein budget with the 
transport restricted to the admissible displacements. Either budget
together with the low-rank perturbation model \eqref{eq:intro-model} defines a structured ambiguity set
$\mathcal A_{U,V}(\PP)$ of distributions of $(X+UCV^\top,Y)$. We fix $\theta$
throughout and study the inner problem
\[
 \sup_{\QQ\in\mathcal A_{U,V}(\PP)}\E_{\xi'\sim \QQ}[\ell(\theta;\xi')],
\]
namely its value, which we call the structured worst-case risk, and the distributions that nearly attain it. These distributions are worst-case scenarios for stress
testing $\theta$. Two examples illustrate settings in which such structured perturbations arise naturally.

\begin{example}[Power system data]\label{ex:power}
Let $X\in\R^{m_1\times m_2}$ record regional net load over a day, with
regions indexed by rows and time periods by columns. Multi-site load data
often exhibit recurring spatial and temporal patterns that motivate low-rank
representations \citep{mateos2013load}. Low-dimensional structure has also
been observed in measurements from time-synchronized grid sensors
\citep{xie2014pmu,gao2016pmu}. We represent coherent load perturbations by
shared spatial patterns $U$ and temporal profiles $V$. For example, a demand
shift with spatial weights $u$ and daily profile $v$ has the rank-one form
$Z=uv^\top$, while multiple such modes give $Z=UCV^\top$. The factors
$U$ and $V$ capture patterns shared across observations, whereas the core
distribution may vary with the nominal observation. A separable budget
$B=B_{c}\otimes B_{r}\succ0$ can further encode temporal and spatial
second-moment geometry, consistent with separable covariance models
\citep{genton2007separable}.
\end{example}

\begin{example}[Image collections]\label{ex:images}
Let $X\in\R^{m_1\times m_2}$ be an image from a collection acquired on a
common grid, and let $\theta$ be a fixed classifier or reconstruction model.
Image collections often exhibit low-dimensional 
variation across samples. For example, images of a fixed convex Lambertian
object under varying distant illumination lie close to a nine-dimensional
linear subspace \citep{basri2003lambertian}. Matrix-aware dimension reduction
further exploits structure shared across a collection: two-dimensional
principal component analysis uses one-sided projections
\citep{yang2004twodpca}, while generalized low-rank approximation learns
common row and column subspaces \citep{ye2005glram}. 
A perturbation $Z=UCV^\top$ then moves every image along directions shared by the collection, by an amount $C$ that can differ from image to image.
\end{example}

The shared support restricts the perturbations available to the adversary and therefore changes the ambiguity set itself. This leads to two distinct questions. First, when does the support restriction preserve the unrestricted worst-case risk, and how large is the resulting risk loss otherwise? Second, for a fixed support, how can one construct a feasible distribution whose risk is close to the corresponding structured worst-case risk? 
To address these questions, we first reduce both budgets to the core coordinates and establish strong duality together with a characterization of optimal distributions for the structured ambiguity sets (Section~\ref{sec:model}). Building on this formulation, our main contributions are as follows.

\begin{itemize}[leftmargin=*,itemsep=2pt]
\item \emph{Excluded risk and lossless supports.}
Unlike principal components, which keep the directions of dominant nominal
variation, the supports that preserve worst-case risk are determined by the
loss and the budget.
Under the moment budget, Theorem~\ref{thm:ridge} gives
two-sided bounds on the excluded risk in terms of the feature budget lost by
the support, and shows that the restriction is lossless when the support
contains the budget-transformed feature directions, with this condition also
being necessary for strongly convex feature losses.
Under the Wasserstein budget, Theorem~\ref{prop:energy} gives analogous
two-sided bounds in terms of the energy of the full-space worst-case
responses missed by the support, with losslessness when all such
responses lie in the prescribed support. Corollaries~\ref{cor:least}
and~\ref{cor:w-least} identify the minimum lossless row and column
dimensions.

\item \emph{Worst-case generation.}
For a fixed support, a quadratic curvature correction turns the deterministic
response problem into a concave program over the original budget set. Its
optimal value upper bounds the structured worst-case risk and 
equals the dual value restricted to multipliers of at least the curvature correction (Proposition~\ref{prop:exact-shift-recovery}).
Projected gradient ascent keeps
every iterate feasible, globally converges at rate $O(1/n)$ in both objective value and
first-order measure, and returns at every iteration a computable interval
containing the structured worst-case risk
(Theorem~\ref{thm:generator-first-order} and
Proposition~\ref{prop:first-order-certificate}). 
This interval closes exactly when the correction vanishes in the limit, and the iterates converge linearly when the budget constraint is sufficiently active, in the sense that the optimal multiplier exceeds the curvature correction (Corollary~\ref{cor:strict-generation}).
\end{itemize}

\subsection{Related work}

\emph{Structured uncertainty and perturbation geometry.}
Moment-based and transportation-based DRO provide the duality and
worst-case distribution theory underlying our two budgets
\citep{delage2010,xu2018,mohajerin2018wasserstein,
blanchetMurthy2019,gaoKleywegt2023}.
More broadly, structure has been used to specify which perturbations
are admissible or costly. Budgeted robust optimization controls
aggregate parameter deviations \citep{bertsimasSim2004}, while
robust factor models account for uncertainty in factor-model
parameters \citep{goldfarbIyengar2003}.
Structured DRO for inverse problems restricts perturbations to
reflect the data-acquisition process \citep{van2026distributionally}.
Low-dimensional or anisotropic transport geometry also appears in
representation-aware costs \citep{wang2025representation} and
data-adaptive robust principal component analysis
\citep{xu2026distributionally}.
Structured adversarial attacks impose matrix structure on individual
perturbations through low-rank factorizations
\citep{savostianova2024lowrank} or nuclear-norm distortion sets
\citep{kazemi2021structured}.
Closer to our setting, \citet{dohmatob2023origins} study
sample-dependent perturbations within a common low-dimensional
subspace and relate fooling rates to margins and gradient alignment.
Shared row and column structure also appears in models of
matrix-valued data, including matrix factor models with common
loadings and a time-varying core \citep{wang2019factor}, and
CoreFlow, which learns shared row and column subspaces while
modeling sample-specific variation through a low-dimensional core
distribution \citep{wu2026coreflow}.
Here, we quantify the loss in worst-case risk caused by restricting
displacements to a prescribed shared row--column support.
With the nominal distribution, loss, and budget fixed, we compare
the unrestricted and support-restricted worst-case risks,
characterize when they coincide, and bound their gap otherwise.

\emph{Dimension reduction and risk preservation.}
Principal component analysis (PCA) based approaches construct reduced-dimensional approximations
for DRO under moment-based and combined moment--Wasserstein ambiguity
sets \citep{cheng2018pca,cheramin2022approximations}.
\citet{jiang2026optimized} optimize the reduced representation jointly
with the moment DRO problem, obtaining low-dimensional semidefinite
approximations that can attain the original optimal value.
On the Wasserstein side, \citet{paty2019subspace} introduce projection-
and subspace-robust Wasserstein distances based on low-dimensional
projections. For losses depending on low-dimensional affine features,
\citet{nietert2024robust} establish an exact feature-space reduction
of the Wasserstein worst-case expectation and use max-sliced
Wasserstein control to derive statistical guarantees.
In statistics, sufficient dimension reduction seeks representations
of predictors that preserve the conditional distribution of the
response \citep{cook2007dimension}. For matrix-valued predictors,
dimension folding \citep{li2010dimension} uses bilinear reductions
$U^\top XV$, a row--column structure closely related to our
perturbation model \eqref{eq:intro-model}.
We ask a related preservation question for worst-case risk:
which shared row and column spaces preserve the unrestricted
worst-case value when admissible displacements are restricted to them?
With the nominal distribution, decision, and budget fixed, we
characterize lossless supports and quantify the risk excluded by
a prescribed support. In our setting, losslessness depends jointly
on the loss and the uncertainty budget.

\emph{Worst-case generation.}
Worst-case distributions have been constructed through penalized adversarial
responses, iterative schemes, and learned transport maps. \citet{sinha2018}
solve pointwise Wasserstein-penalized response problems and use the
resulting perturbations for distributionally robust training.
\citet{blanchet2022optimal} exploit the structure of the worst-case adversarial model under affine decision rules to solve optimal-transport DRO by stochastic gradient methods at the cost of the non-robust problem.
\citet{xu2024flow,zhu2024iterative,cheng2025} represent
continuous worst-case distributions through transport maps and optimize
these maps directly. We consider a different setting in which the decision and a shared
row--column support are fixed, and the adversary is constrained by a
prescribed moment or Wasserstein budget. Our goal is to generate
feasible perturbations within this support while certifying their
suboptimality relative to the original conditional distribution problem.

\subsection{Organization}
Section~\ref{sec:model} defines the structured ambiguity sets and
establishes duality for the associated worst-case risk problems.
Sections~\ref{sec:matrix-structure} and~\ref{sec:w-structure}
quantify the excluded risk and characterize lossless supports under
the moment and Wasserstein budgets, respectively.
Section~\ref{sec:recovery} develops a worst-case generation algorithm
with computable risk certificates for a fixed support, and
Section~\ref{sec:experiments} reports experiments on spatial
precipitation forecasting.

\subsection{Notation}

Our notation is mostly standard. For matrices, $\langle A,B\rangle_F:=\tr(A^\top B)$ denotes the Frobenius
inner product and $\|\cdot\|_F$ the associated norm, and we omit the
subscript when no confusion arises. Vector norms are Euclidean, and
$\|\cdot\|_{\rm op}$ denotes the induced operator norm.
We write $\mathbb S_{++}^k$ for the set of $k\times k$ symmetric
positive definite matrices. For symmetric matrices $A$ and $B$,
$A\succeq B$ (resp. $A\succ B$) means that $A-B$ is positive semidefinite (resp. positive definite).
The symbols $\lambda_{\min}(A)$
and $\lambda_{\max}(A)$ denote the smallest and largest eigenvalues of a
symmetric matrix $A$, and for $A\succeq0$, $A^{1/2}$ denotes its unique
positive semidefinite square root. The symbol $\dagger$ denotes the
Moore--Penrose pseudoinverse. The operator $\vecop$ stacks the columns of
a matrix, $\matop$ denotes its inverse with the matrix dimensions
understood from context, and $\otimes$ denotes the Kronecker product. In
particular, $\vecop(UCV^\top)=(V\otimes U)\vecop(C)$. We write $L^2(\PP)$
for the space of square-integrable real functions of the nominal
observation $\xi$. For a square-integrable random core matrix $C$ with conditional
distribution $q(\cdot\mid\xi)$, we write
$\|C\|_{L^2}^2:=\E_q\|C\|_F^2$. Finally, $\delta_x$ denotes the point mass at
$x$.

\section{Low-Rank Ambiguity Sets}\label{sec:model}

Let $\xi=(X,Y)$ have nominal distribution $\PP$ on
$\R^{m_1\times m_2}\times\mathcal Y$, where $\mathcal Y$ is a standard Borel
space, and at a fixed decision $\theta$ write
\[
 f_\xi(Z):=\ell(\theta;X+Z,Y).
\]
In the inner maximization problem under the perturbation model \eqref{eq:intro-model}, the optimization variable is a conditional distribution $q(dC\mid\xi)$ of the core matrix $C$ in $Z=UCV^\top$. We require $\E_q\|C\|_F^2<\infty$, so that both budgets in
\eqref{eq:intro-budgets} are well defined.
Write
\begin{equation}\label{eq:core-coordinates}
 S:=V\otimes U,\quad k:=rs,\quad d:=m_1m_2,
 \quad\vecop(UCV^\top)=S\vecop(C).
\end{equation}
Then $S^\top S=I_k$ and $\|UCV^\top\|_F=\|C\|_F$. We first have the following lemma.

\begin{lemma}[Induced core budget]\label{prop:core-budget}
Let
\begin{equation}\label{eq:induced-budget}
 B_S:=(S^\top B^{-1}S)^{-1}.
\end{equation}
Then, for every symmetric $M\in\R^{k\times k}$,
\begin{equation}\label{eq:core-budget-equivalence}
 SMS^\top\preceq B\quad\Longleftrightarrow\quad M\preceq B_S.
\end{equation}
In particular, $SB_SS^\top\preceq B$, and $B_S$ is the largest $M$ in the
Loewner order with $SMS^\top\preceq B$.
\end{lemma}
\begin{proof}[Proof of Lemma \ref{prop:core-budget}]
First, the matrix $S^\top B^{-1}S$ is positive definite because $B\succ0$ and
$S$ has full column rank. Thus $B_S$ in \eqref{eq:induced-budget} is
well defined. Define
$
 \tilde{S}:=B^{-1/2}SB_S^{1/2}$ and $\tilde{M}:=B_S^{-1/2}MB_S^{-1/2}$.
By the definition of $B_S$,
$$
\tilde{S}^\top \tilde{S}
=
B_S^{1/2}S^\top B^{-1}SB_S^{1/2}
=
I_k,\qquad
\tilde{S}\tilde{M}\tilde{S}^\top
=
B^{-1/2}SMS^\top B^{-1/2}.
$$
Since congruence by the invertible matrix $B^{-1/2}$ preserves the
Loewner order,
$$
SMS^\top\preceq B
\quad\Longleftrightarrow\quad
B^{-1/2}SMS^\top B^{-1/2}\preceq I_d
\quad\Longleftrightarrow\quad
\tilde{S}\tilde{M}\tilde{S}^\top\preceq I_d.
$$
It remains to show that
$$
\tilde{S}\tilde{M}\tilde{S}^\top\preceq I_d
\quad\Longleftrightarrow\quad
\tilde{M}\preceq I_k.
$$
If $\tilde{S}\tilde{M}\tilde{S}^\top\preceq I_d$, multiplying by $\tilde{S}^\top$ and $\tilde{S}$ and using
$\tilde{S}^\top \tilde{S}=I_k$ gives $\tilde{M}\preceq I_k$. Conversely, if $\tilde{M}\preceq I_k$,
then
$
I_d-\tilde{S}\tilde{M}\tilde{S}^\top
=
(I_d-\tilde{S}\tilde{S}^\top)+\tilde{S}(I_k-\tilde{M})\tilde{S}^\top
\succeq0,
$
because $\tilde{S}\tilde{S}^\top$ is an orthogonal projector. 
Finally,
$\tilde{M}\preceq I_k$ is equivalent to $M\preceq B_S$ by congruence with
$B_S^{1/2}$. This proves \eqref{eq:core-budget-equivalence}.
Taking $M=B_S$ gives $SB_SS^\top\preceq B$.
\end{proof}

To ensure that the worst-case expectations are finite under the two
second-moment budgets, we impose a mild quadratic growth condition on the
loss. It allows the first-order growth coefficient to depend on the nominal
observation, while the quadratic term is uniform across observations.

\begin{assumption}[Loss regularity and growth]\label{ass:growth}
The map $(\xi,Z)\mapsto f_\xi(Z)$ is jointly measurable, and
$Z\mapsto f_\xi(Z)$ is continuous for $\PP$-almost every $\xi$.
Moreover, there exist a nonnegative $a\in L^2(\PP)$ and a constant
$L\ge0$ such that
\begin{equation}\label{eq:growth}
 |f_\xi(Z)-f_\xi(0)|
 \le
 a(\xi)\|Z\|_F+\frac{L}{2}\|Z\|_F^2
\end{equation}
for $\PP$-almost every $\xi$ and every $Z$, with
$\E|f_\xi(0)|<\infty$.
\end{assumption}

For a conditional core distribution $q$ in \eqref{eq:intro-model}, define its
expected loss and core second moment by
\begin{equation}\label{eq:core-loss-moment}
 J_S(q):=\E_q f_\xi(UCV^\top),\qquad
 M(q):=\E_q[\vecop(C)\vecop(C)^\top].
\end{equation}
Under Assumption~\ref{ass:growth}, the Cauchy--Schwarz inequality shows that
$J_S(q)$ is finite whenever $\E_q\|C\|_F^2<\infty$.
By \eqref{eq:core-coordinates},
$\E_q[\vecop(Z)\vecop(Z)^\top]=SM(q)S^\top$, so Lemma~\ref{prop:core-budget}
shows that the moment budget in \eqref{eq:intro-budgets} is equivalent to
$M(q)\preceq B_S$. Since $\|UCV^\top\|_F=\|C\|_F$, the Wasserstein budget is
equivalent to $\E_q\|C\|_F^2\le\rho^2$. Then the two structured worst-case risks are
\begin{align}
 R_S^{\rm M}(B)&:=\sup_{q:\,M(q)\preceq B_S}J_S(q),\label{eq:matrix-risk}\\
 R_S^{\rm W}(\rho)&:=\sup_{q:\,\E_q\|C\|_F^2\le\rho^2}J_S(q).\label{eq:wasserstein-risk}
\end{align}
Taking $U=I_{m_1}$ and $V=I_{m_2}$ gives $S=I_d$ and $B_S=B$, recovering the
full-space problems, and we denote the corresponding risks by $R^{\rm M}(B)$
and $R^{\rm W}(\rho)$. When the statement applies to
both moment and Wasserstein cases, we write simply $R$ and $R_S$.
Every structured feasible distribution is
feasible for the full-space problem with the same budget, so $R_S\le R$, and
for every feasible $q$,
\begin{equation}\label{eq:total-error}
 R-J_S(q)=\Delta_S+(R_S-J_S(q)),\qquad \Delta_S:=R-R_S,
\end{equation}
where the excluded risk $\Delta_S$ and the generation error $R_S-J_S(q)$ are
both nonnegative.

\subsection{Duality and optimality conditions}\label{sec:duality}
We next derive dual representations for the two structured worst-case
problems. The moment budget is priced by a positive semidefinite $k\times k$
matrix multiplier $\Lambda$, and the Wasserstein budget by a nonnegative scalar
multiplier $\lambda$. These multipliers 
will
also be used to characterize optimal conditional distributions. For a nominal
observation $\xi$, define the priced losses
\[
 \psi_{\xi,\Lambda}(C):=f_\xi(UCV^\top)-\vecop(C)^\top\Lambda\vecop(C),
 \qquad
 \psi_{\xi,\lambda}(C):=f_\xi(UCV^\top)-\lambda\|C\|_F^2,
\]
and write $\psi_\xi$ when the budget and multiplier are clear from context.
The dual functions are
\begin{align}
 \mathcal D_S^{\rm M}(\Lambda)&:=\tr(\Lambda B_S)+\E\sup_C\psi_{\xi,\Lambda}(C),
 \label{eq:dual-value}\\
 \mathcal D_S^{\rm W}(\lambda)&:=\lambda\rho^2+\E\sup_C\psi_{\xi,\lambda}(C).
 \label{eq:dual-value-w}
\end{align}
Since $\sup_C\psi_\xi(C)\ge f_\xi(0)$ and, by continuity, the supremum can be
taken over a countable dense set, both dual functions are well defined.
To characterize optimality, we introduce two nonnegative terms in the primal--dual gap:
the unused budget weighted by its
multiplier (the priced budget slack) and the loss from not selecting a
pointwise optimal priced response (priced response error),
\begin{align}
 s(q)&:=
 \begin{cases}
 \tr\bigl[\Lambda\bigl(B_S-M(q)\bigr)\bigr] & \text{(moment budget)},\\
 \lambda\bigl(\rho^2-\E_q\|C\|_F^2\bigr) & \text{(Wasserstein budget)},
 \end{cases}
 \label{eq:budget-slack}\\
 e(q)&:=\E_q\Bigl[\sup_{C'}\psi_\xi(C')-\psi_\xi(C)\Bigr].
 \label{eq:response-error}
\end{align}

The following proposition establishes strong duality and shows that
the primal--dual gap is the sum of the priced budget slack and
priced response error.

\begin{proposition}[Strong duality and optimality characterization]
\label{prop:duality}
Under Assumption~\ref{ass:growth},
\begin{equation}\label{eq:strong-duality}
 R_S^{\rm M}=\min_{\Lambda\succeq0}\mathcal D_S^{\rm M}(\Lambda),\qquad
 R_S^{\rm W}=\min_{\lambda\ge0}\mathcal D_S^{\rm W}(\lambda),
\end{equation}
with both dual minima attained.
Moreover, for any feasible conditional distribution \(q\) and any
multiplier with finite dual value,
\begin{equation}\label{eq:gap-identity}
 \mathcal D_S^{\rm M}(\Lambda)-J_S(q)=s(q)+e(q),\qquad \mathcal D_S^{\rm W}(\lambda)-J_S(q)=s(q)+e(q).
\end{equation}
Consequently, at any dual optimum,
a feasible $q$ is worst-case optimal if and only if $s(q)=e(q)=0$.
\end{proposition}

\begin{proof}[Proof of Proposition \ref{prop:duality}]
Throughout this proof, let
\[
 \mathcal V(D):=\sup_{M(q)\preceq D}J_S(q),\qquad D\in\mathbb S^k_{++}.
\]
By the definitions in
\eqref{eq:core-loss-moment}, and the growth estimate \eqref{eq:growth}, every
finite second moment $q$ satisfies
\[
 J_S(q)-\E f_\xi(0)
 \le\E_q[a(\xi)\|C\|_F]+\frac L2\E_q\|C\|_F^2
 \le\|a\|_{L^2(\PP)}\sqrt{\tr M(q)}+\frac L2\tr M(q).
\]
Since the zero displacement is feasible and $\tr M(q)\le\tr D$ whenever
$M(q)\preceq D$,
\begin{equation}\label{eq:value-finite}
 \E f_\xi(0)\le\mathcal V(D)\le\E f_\xi(0)
 +\|a\|_{L^2(\PP)}\sqrt{\tr D}+\frac L2\tr D<\infty.
\end{equation}
For $D_1,D_2\succ0$, $0<t<1$, and $\varepsilon>0$, choose distributions $q_i$, $i=1,2$, 
such that $M(q_i)\preceq D_i$ and $J_S(q_i)\ge\mathcal V(D_i)-\varepsilon$.
Their conditional mixture $q_t=t q_1+(1-t)q_2$ 
satisfies
\[
 M(q_t)=t M(q_1)+(1-t)M(q_2),\qquad
 J_S(q_t)=t J_S(q_1)+(1-t)J_S(q_2).
\]
It follows that
\[
 \mathcal V(t D_1+(1-t)D_2)\ge J_S(q_t)
 \ge t\mathcal V(D_1)+(1-t)\mathcal V(D_2)-\varepsilon.
\]
Letting $\varepsilon\downarrow0$ proves concavity of $\mathcal V$. Inclusion of the
feasible sets also shows that $\mathcal V(D_1)\le\mathcal V(D_2)$ whenever
$0\prec D_1\preceq D_2$.

Next we construct a supporting multiplier and show that it is positive
semidefinite. Extend $-\mathcal V$ by $+\infty$ outside $\mathbb S^k_{++}$. By
\eqref{eq:value-finite} and the preceding concavity argument, it is a
proper convex function whose effective domain is open. Since
$B_S\succ0$ by Lemma~\ref{prop:core-budget},
Rockafellar's subgradient theorem \citep[Theorem~23.4]{rockafellar1970} gives a
subgradient $-\Lambda_*\in\partial(-\mathcal V)(B_S)$, and
\begin{equation}\label{eq:value-support}
 \mathcal V(D)\le\mathcal V(B_S)+\tr[\Lambda_*(D-B_S)],\quad D\succ0.
\end{equation}
For any \(u\in\R^k\) and \(t'>0\), we have
\(B_S\preceq B_S+t'uu^\top\). By monotonicity of \(\mathcal V\) and
\eqref{eq:value-support},
\[
 0\le\mathcal V(B_S+t'uu^\top)-\mathcal V(B_S)\le t'u^\top\Lambda_*u.
\]
Since \(t'>0\) and \(u\) is arbitrary, \(\Lambda_*\succeq0\).

Now take any finite second moment conditional distribution $q$.
For every $\delta>0$, $q$ is feasible at
$D=M(q)+\delta I_k\succ0$. Then
\eqref{eq:value-support} gives 
\[
 J_S(q)\le\mathcal V(M(q)+\delta I_k)
 \le\mathcal V(B_S)+\tr[\Lambda_*(M(q)+\delta I_k-B_S)].
\]
After rearranging and letting $\delta\downarrow0$, we obtain
\begin{equation}\label{eq:lagrange-upper}
 J_S(q)-\tr[\Lambda_*M(q)]
 \le\mathcal V(B_S)-\tr(\Lambda_*B_S).
\end{equation}
To pass from the pointwise supremum to feasible responses, fix
$\Lambda\succeq0$. By Assumption~\ref{ass:growth},
$(\xi,C)\mapsto\psi_{\xi,\Lambda}(C)$ is jointly measurable and, after a
$\PP$-null modification, continuous in $C$, so $-\psi_{\xi,\Lambda}$ is a
normal integrand with $\E|\psi_{\xi,\Lambda}(0)|<\infty$. The
interchange theorem \citep[Theorem~14.60]{rockafellar2009variational}
 gives
\begin{equation}\label{eq:finite-response-envelope}
 \E\sup_C\psi_{\xi,\Lambda}(C)
 =\sup_T\E\,\psi_{\xi,\Lambda}(T(\xi)),
\end{equation}
where the supremum is over measurable maps
$T:\xi\mapsto T(\xi)\in\R^{r\times s}$ with $\E\|T(\xi)\|_F^2<\infty$,
which form a decomposable space. Each such $T$ induces the deterministic
conditional distribution $q_T(dC\mid\xi)=\delta_{T(\xi)}(dC)$ with finite
second moment, so \eqref{eq:lagrange-upper} applies to $q_T$ at
$\Lambda=\Lambda_*$. Taking the supremum over $T$ and using
\eqref{eq:finite-response-envelope},
\[
 \E\sup_C\psi_{\xi,\Lambda_*}(C)\le\mathcal V(B_S)-\tr(\Lambda_*B_S),
\]
so $\mathcal D_S^{\rm M}(\Lambda_*)
\le
\mathcal V(B_S)$.
For the reverse inequality, every feasible distribution and
every $\Lambda\succeq0$ satisfy
\[
 J_S(q)=\E_q\psi_{\xi,\Lambda}(C)+\tr[\Lambda M(q)]
 \le\E\sup_C\psi_{\xi,\Lambda}(C)+\tr(\Lambda B_S)=\mathcal D_S^{\rm M}(\Lambda),
\]
where the inequality uses $M(q)\preceq B_S$ and
$\Lambda\succeq0$.
Taking the supremum over the distributions in
\eqref{eq:matrix-risk}, whose value is $R_S^{\rm M}=\mathcal V(B_S)$, and then
the infimum over multipliers gives
\[
 \mathcal V(B_S)\le\inf_{\Lambda\succeq0}\mathcal D_S^{\rm M}(\Lambda)
 \le \mathcal D_S^{\rm M}(\Lambda_*)\le\mathcal V(B_S).
\]
This proves matrix strong duality and multiplier attainment in
\eqref{eq:strong-duality}.

For the scalar budget, define
$
 \mathcal W(t):=\sup_{\E_q\|C\|_F^2\le t}J_S(q)$, $t>0$.
The similar growth and mixture arguments above give
\[
 \E f_\xi(0)\le\mathcal W(t)\le\E f_\xi(0)+\|a\|_{L^2(\PP)}\sqrt t+Lt/2,
\]
and show that $\mathcal W$ is concave and nondecreasing. Its supporting slope
at the interior point $t=\rho^2$ can therefore be chosen as $\lambda_*\ge0$.
Applying the scalar supporting inequality at
$t=\E_q\|C\|_F^2+\delta$ and then letting $\delta\downarrow0$ yields
$
 J_S(q)-\lambda_*\E_q\|C\|_F^2\le\mathcal W(\rho^2)-\lambda_*\rho^2$
for every finite moment kernel. The construction leading to
\eqref{eq:finite-response-envelope} also applies to $\psi_{\xi,\lambda_*}$ and gives
\[
 \E\sup_C\psi_{\xi,\lambda_*}(C)\le\mathcal W(\rho^2)-\lambda_*\rho^2.
\]
Weak duality for \eqref{eq:wasserstein-risk} and \eqref{eq:dual-value-w}
closes the same chain, since $\mathcal W(\rho^2)=R_S^{\rm W}$:
\[
 \mathcal W(\rho^2)\le\inf_{\lambda\ge0}\mathcal D_S^{\rm W}(\lambda)
 \le \mathcal D_S^{\rm W}(\lambda_*)\le\mathcal W(\rho^2).
\]
This proves the scalar part of \eqref{eq:strong-duality} and its
multiplier attainment.

Finally, substitute \eqref{eq:budget-slack} and \eqref{eq:response-error}
into the difference between \eqref{eq:dual-value} and
\eqref{eq:core-loss-moment}. For a matrix price this gives
\[
 s(q)+e(q)
 =\tr(\Lambda B_S)-\tr[\Lambda M(q)]
 +\E\sup_{C'}\psi_\xi(C')-\E_q\psi_\xi(C)=\mathcal D_S^{\rm M}(\Lambda)-J_S(q).
\]
Replacing the trace terms by $\lambda\rho^2$ and
$\lambda\E_q\|C\|_F^2$ gives the scalar identity. 
Thus \eqref{eq:gap-identity} holds. At a dual optimal
multiplier, \eqref{eq:strong-duality} identifies the left-hand side of
\eqref{eq:gap-identity} with $R_S-J_S(q)$, so a feasible $q$ is optimal
exactly when both terms vanish.
\end{proof}

Proposition~\ref{prop:duality} parallels the familiar strong duality
structure for moment-based \citep{xu2018} and optimal-transport
\citep{blanchetMurthy2019,gaoKleywegt2023} ambiguity sets, here with
displacements restricted to a prescribed support. At an optimal multiplier,
the gap identity shows that a worst-case distribution must be supported on
pointwise maximizers of the priced loss and have zero priced budget slack.
A pointwise maximizing deterministic response need not satisfy the budget
complementarity condition $s(q)=0$. In such cases, randomization among pointwise
maximizers may be needed to attain the primal optimum.

From now on, we take the nominal distribution to be finitely supported, with observations $\xi_1,\dots,\xi_N$ and masses $p_i>0$, $\sum_ip_i=1$, as for an empirical distribution, and we write $f_i:=f_{\xi_i}$. 
The conditional laws $q(\cdot\mid\xi_i)$ are not restricted to
point masses and may remain randomized.

\section{Moment Budgets}\label{sec:matrix-structure}

In this section, we study how a fixed shared support affects the worst-case risk under the
moment budget. To obtain an exact characterization, suppose that the loss
depends on the displacement through fixed row and column features. Let
\(U_0\in\R^{m_1\times r_0}\) and \(V_0\in\R^{m_2\times s_0}\) have orthonormal
columns, and set $S_0:=V_0\otimes U_0$ and $k_0:=r_0s_0$. 
For twice
continuously differentiable \(h_i:\R^{r_0\times s_0}\to\R\), assume
\begin{equation}\label{eq:ridge}
 f_i(Z)=h_i(U_0^\top ZV_0)
 \quad\text{and}\quad
 -LI_{k_0}\preceq\nabla^2h_i(C_0)\preceq LI_{k_0}
 \ \text{for all } i \text{ and } C_0\in\R^{r_0\times s_0}.
\end{equation}
Thus the loss sees the displacement only through its coordinates $S_0^\top\vecop(Z)$ in $\range(S_0)$.
Write
$g_0:=(\sum_i p_i\|\nabla h_i(0)\|_F^2)^{1/2}$
for the root-mean-square feature gradient at zero displacement. For a feature
budget $D'\succeq0$, define
\begin{equation}
\label{eq:feature-value-definition}
 \mathcal V_0(D'):=
 \sup_{\E[\vecop(C_0)\vecop(C_0)^\top]\preceq D'}
 \sum_i p_i\,\E[h_i(C_0)\mid i],
\end{equation}
where the supremum is over conditional distributions of
$C_0\in\R^{r_0\times s_0}$ given the nominal index $i$, and the expectation
in the constraint includes both $i\sim(p_1,\dots,p_N)$ and the conditional law of
$C_0$. 
Problem \eqref{eq:feature-value-definition} is the full-space moment problem of
Section~\ref{sec:model} in dimension $k_0$ with losses $h_i$, so
Proposition~\ref{prop:duality} applies to it. Let $D:=S_0^\top BS_0$,
which is positive definite since $B\succ0$ and $S_0$ has orthonormal
columns, and let $\Lambda_0\succeq0$ be a dual-optimal multiplier for
$\mathcal V_0(D)$, that is,
\begin{equation}\label{eq:feature-dual}
\mathcal V_0(D)=\tr(\Lambda_0D)
+\sum_ip_i\sup_{C_0}\bigl\{h_i(C_0)-\vecop(C_0)^\top\Lambda_0\vecop(C_0)\bigr\}.
\end{equation}

The main theorem characterizes how the prescribed support reduces the
effective feature budget and quantifies the resulting loss in worst-case risk.

\begin{theorem}[Excluded risk under moment budgets]
\label{thm:ridge}
Let $B\succ0$ and assume \eqref{eq:ridge}. For an orthonormal support
$S$, let
\begin{equation*}
D:=S_0^\top BS_0,
\qquad
D_S:=S_0^\top S B_S S^\top S_0.
\end{equation*}
Then $0\preceq D_S\preceq D$, and the full-space and restricted risks
admit the exact representations
\[
R^{\rm M}(B)=\mathcal V_0(D),
\qquad
R_S^{\rm M}(B)=\mathcal V_0(D_S),
\]
and the excluded risk satisfies
\begin{equation}\label{eq:ridge-upper}
\tr[\Lambda_0(D-D_S)]\le\Delta_S^{\rm M}\le c_D\tr(D-D_S),
\end{equation}
where $c_D:=(g_0+L\sqrt{\tr D})/\sqrt{\lambda_{\min}(D)}$. 
If, in addition, $\nabla^2h_i(C_0)\succeq\alpha I_{k_0}$ for all $i$ and
$C_0$ and some $\alpha\ge0$, then $\Lambda_0\succeq\frac\alpha2I_{k_0}$, and
hence
\begin{equation}\label{eq:ridge-bound}
\Delta_S^{\rm M}\ge\frac\alpha2\tr(D-D_S).
\end{equation}
The inclusion $\range(BS_0)\subseteq\range(S)$ is sufficient for
$\Delta_S^{\rm M}=0$, and is also necessary when $\Lambda_0\succ0$, in
particular when $\alpha>0$.
\end{theorem}

\begin{proof}[Proof of Theorem \ref{thm:ridge}]
For the proof, we identify each feature matrix $C_0$ with its vectorization
$c_0=\vecop(C_0)$ and write $h_i(c_0):=h_i(\matop(c_0))$.
Since $B\succ0$ and $S_0$ has orthonormal columns, we know $D\succ0$.

We first show that the support restriction is exactly represented in the
feature space by replacing the full budget $D$ with $D_S$. Define
$\tilde{S}:=B^{-1/2}SB_S^{1/2}$ and
$F:=S_0^\top SB_S^{1/2}$.
Then $\tilde{S}^\top \tilde{S}=I_k$ and $D_S=FF^\top\succeq0$. Moreover,
\begin{equation}\label{eq:ridge-budget-deficit}
D-D_S
=S_0^\top B^{1/2}(I-\tilde{S}\tilde{S}^\top)B^{1/2}S_0
\succeq0.
\end{equation}
By Lemma~\ref{prop:core-budget}, a conditional distribution of the core
vector $w$ satisfies the ambient moment constraint if and only if
$\E[ww^\top]\preceq B_S$. Its feature vector $c_0=S_0^\top Sw$ therefore
satisfies
\[
\E[c_0c_0^\top]
\preceq S_0^\top SB_SS^\top S_0=D_S.
\]
Conversely, suppose $\E[c_0c_0^\top]\preceq D_S=FF^\top$. For every
$v\in\ker F^\top$,
$
\E(v^\top c_0)^2\le v^\top FF^\top v=0$.
Applying this to a finite basis of $\ker F^\top$ gives
$c_0\in\range F$ almost surely. Define
$w=B_S^{1/2}F^\dagger c_0$.
Since $FF^\dagger$ and $F^\dagger F$ are the orthogonal projectors onto
$\range F$ and $\range F^\top$, respectively,
\[
S_0^\top Sw=FF^\dagger c_0=c_0\quad\text{and}\quad
\E[ww^\top]\preceq
B_S^{1/2}F^\dagger FF^\top(F^\dagger)^\top B_S^{1/2}=
B_S^{1/2}(F^\dagger F)B_S^{1/2}
\preceq B_S.
\]
Thus every $D_S$-feasible feature distribution is the image of a
$B_S$-feasible core distribution under $w\mapsto S_0^\top Sw$, and by
\eqref{eq:ridge} the two have the same expected loss. Hence the core problem \eqref{eq:matrix-risk} and the
feature problem with budget $D_S$ take the same objective values, so
$R_S^{\rm M}(B)=\mathcal V_0(D_S)$. Taking $S=I_d$, for which $B_S=B$ and
$D_S=D$, gives the full-space identity. Hence
\begin{equation}\label{eq:feature-values}
R^{\rm M}(B)=\mathcal V_0(D),
\qquad
R_S^{\rm M}(B)=\mathcal V_0(D_S),
\end{equation}
and $\Delta_S^{\rm M}\ge0$ because $D_S\preceq D$ and $\mathcal V_0$ is
monotone.

The deficit $D-D_S$ in \eqref{eq:ridge-budget-deficit} is a Gram matrix, so
\begin{equation}\label{eq:feature-budget-exact}
\begin{aligned}
D=D_S
&\quad\Longleftrightarrow\quad
(I-\tilde{S}\tilde{S}^\top)B^{1/2}S_0=0\\
&\quad\Longleftrightarrow\quad
\range(B^{1/2}S_0)\subseteq\range(B^{-1/2}S)\\
&\quad\Longleftrightarrow\quad
\range(BS_0)\subseteq\range(S).
\end{aligned}
\end{equation}
Here $\range \tilde{S}=\range(B^{-1/2}S)$, and the last equivalence uses
invertibility of $B^{1/2}$. Together with \eqref{eq:feature-values},
this proves that the range inclusion is sufficient for
$\Delta_S^{\rm M}=0$.

Now we handle the upper estimate. Take any conditional feature distribution with
$M=\E[c_0c_0^\top]\preceq D$, and let
$
y:=D_S^{1/2}D^{-1/2}c_0$,
$H:=D^{1/2}-D_S^{1/2}$.
Then
\begin{equation*}
\begin{aligned}
\E[yy^\top]
&=
D_S^{1/2}D^{-1/2}MD^{-1/2}D_S^{1/2}
\preceq D_S,\quad \E\|c_0-y\|^2=
\tr(HD^{-1/2}MD^{-1/2}H^\top)
\le \|H\|_F^2.
\end{aligned}
\end{equation*}
Taking the Frobenius inner product of $D-D_S=D^{1/2}H+HD_S^{1/2}$ with $H$ gives
\[
\begin{aligned}
\sqrt{\lambda_{\min}(D)}\,\|H\|_F^2
\le
\tr(H^\top D^{1/2}H)+\tr(H^\top HD_S^{1/2})=
\langle H,D-D_S\rangle
\le \|H\|_F\,\|D-D_S\|_F.
\end{aligned}
\]
Since $D-D_S\succeq0$,
\begin{equation}\label{eq:feature-square-root-bound}
\|H\|_F
\le
\frac{\|D-D_S\|_F}{\sqrt{\lambda_{\min}(D)}}
\le
\frac{\tr(D-D_S)}{\sqrt{\lambda_{\min}(D)}}.
\end{equation}

For $t\in[0,1]$, let $z_t:=(1-t)y+tc_0$. Since
$h_i(c_0)-h_i(y)=\int_0^1\langle\nabla h_i(z_t),c_0-y\rangle\,dt$,
Fubini's theorem and Cauchy--Schwarz give
\[
\E[h_i(c_0)-h_i(y)]
\le
\int_0^1\sqrt{\E\|\nabla h_i(z_t)\|^2}\,dt\,
\sqrt{\E\|c_0-y\|^2}.
\]
Since $\E\|y\|^2\le\tr D_S\le\tr D$ and $\E\|c_0\|^2\le\tr D$, Minkowski's
inequality gives
\[
\sqrt{\E\|z_t\|^2}
\le(1-t)\sqrt{\E\|y\|^2}+t\sqrt{\E\|c_0\|^2}
\le\sqrt{\tr D}.
\]
The gradient estimate
$\|\nabla h_i(z_t)\|\le\|\nabla h_i(0)\|+L\|z_t\|$
and Minkowski's inequality then give
\[
\sqrt{\E\|\nabla h_i(z_t)\|^2}
\le
g_0+L\sqrt{\E\|z_t\|^2}
\le
g_0+L\sqrt{\tr D}\quad
\text{for every}\ t\in[0,1].
\]
Together with $\E\|c_0-y\|^2\le\|H\|_F^2$, this yields
\[
\E[h_i(c_0)-h_i(y)]
\le
(g_0+L\sqrt{\tr D})\|H\|_F.
\]
Combining the preceding
bounds with \eqref{eq:feature-square-root-bound} and 
$
\mathcal V_0(D_S)
\ge
\sum_i p_i\E[h_i(y)\mid i]$,
we obtain
\[
\sum_i p_i\E[h_i(c_0)\mid i]-\mathcal V_0(D_S)
\le
\sum_ip_i\E[h_i(c_0)-h_i(y)\mid i]
\le
c_D\tr(D-D_S).
\]
Taking the supremum over all $D$-feasible distributions and using
\eqref{eq:feature-values} proves the upper bound in \eqref{eq:ridge-upper}.

For the lower estimate, let $c_0$ have any conditional distribution
satisfying $\E[c_0c_0^\top]\preceq D_S$. Since $\Lambda_0\succeq0$ and
$\E[c_0c_0^\top]\preceq D_S$, we have
$\tr(\Lambda_0\E[c_0c_0^\top])\le\tr(\Lambda_0D_S)$, and
\eqref{eq:feature-dual} gives
\[
\begin{aligned}
\sum_i p_i\E[h_i(c_0)\mid i]
&=
\sum_i p_i\E[h_i(c_0)-c_0^\top\Lambda_0c_0\mid i]
+\tr(\Lambda_0\E[c_0c_0^\top])\\
&\le
\sum_i p_i\sup_{c_0'}\{h_i(c_0')-c_0'^\top\Lambda_0c_0'\}
+\tr(\Lambda_0D_S)
=\mathcal V_0(D)-\tr[\Lambda_0(D-D_S)].
\end{aligned}
\]
Since this holds for every $D_S$-feasible conditional distribution of
$c_0$, taking the supremum gives
$\mathcal V_0(D)-\mathcal V_0(D_S)\ge\tr[\Lambda_0(D-D_S)]$, which with
\eqref{eq:feature-values} is the lower bound in \eqref{eq:ridge-upper}.

Now suppose $\nabla^2h_i(c_0)\succeq\alpha I_{k_0}$ for all $i$ and $c_0$
and some $\alpha\ge0$. If $\Lambda_0\not\succeq\frac\alpha2I_{k_0}$, there
is $v$ with $v^\top(\alpha I_{k_0}-2\Lambda_0)v>0$, and then for every $i$,
\[
h_i(tv)-t^2v^\top\Lambda_0v
\ge
h_i(0)+t\langle\nabla h_i(0),v\rangle
+\frac{t^2}{2}\,v^\top(\alpha I_{k_0}-2\Lambda_0)v
\to\infty
\qquad(t\to\infty),
\]
so the right-hand side of \eqref{eq:feature-dual} is infinite, whereas
$\mathcal V_0(D)$ is finite. Hence $\Lambda_0\succeq\frac\alpha2I_{k_0}$.
Since $\Lambda_0-\frac\alpha2I_{k_0}\succeq0$ and $D-D_S\succeq0$, we have
$\tr[\Lambda_0(D-D_S)]\ge\frac\alpha2\tr(D-D_S)$, and
\eqref{eq:ridge-bound} follows from \eqref{eq:ridge-upper}.

Finally, if $\Lambda_0\succ0$ and $\Delta_S^{\rm M}=0$, then
\eqref{eq:ridge-upper} gives $\tr[\Lambda_0(D-D_S)]=0$. Since
$\Lambda_0\succ0$ and $D-D_S\succeq0$, this forces $D=D_S$, and
\eqref{eq:feature-budget-exact} yields $\range(BS_0)\subseteq\range(S)$.
Thus the range inclusion is also necessary when $\Lambda_0\succ0$, in
particular when $\alpha>0$, since then $\Lambda_0\succeq\frac\alpha2I_{k_0}\succ0$.
\end{proof}

\begin{remark}[Loss features and lossless perturbation supports]
\label{rem:loss-features}
The factorization $f_i(Z)=h_i(U_0^\top ZV_0)$ means that each loss is determined by $S_0^\top\vecop(Z)$ alone, analogous
to sufficient representations in dimension reduction and dimension
folding \citep{cook2007dimension,li2010dimension}.
Preserving the worst-case risk is a different matter:
by Theorem~\ref{thm:ridge}, the directions that
matter are those of the budget-transformed space $\range(BS_0)$, which
need not coincide with the loss feature space $\range(S_0)$. The two
coincide exactly when $B\,\range(S_0)\subseteq\range(S_0)$, as in the
isotropic case $B=\sigma^2I_d$ with $\sigma>0$.
\end{remark}

The range condition in Theorem~\ref{thm:ridge} also determines how many common
column and row directions suffice and how to choose their bases. When
$\Lambda_0\succ0$, no smaller shared support is lossless.

\begin{corollary}[Minimal lossless dimensions]
\label{cor:least}
Let $B\succ0$ and assume \eqref{eq:ridge}. 
For $j=1,\ldots,k_0$, let $s_{0j}$ be the $j$th column of $S_0$. Let $U_*\in\R^{m_1\times r_*}$ and $V_*\in\R^{m_2\times s_*}$ have orthonormal columns with
\begin{equation}\label{eq:feature-common-spans}
 \range(U_*)=\sum_{j=1}^{k_0}\range(\matop(Bs_{0j})),
 \qquad
 \range(V_*)=\sum_{j=1}^{k_0}\range(\matop(Bs_{0j})^\top),
\end{equation}
where $r_*$ and $s_*$ are the dimensions of these two subspaces.
Then $S_*:=V_*\otimes U_*$ is lossless, and so is every orthonormal support
$S=V\otimes U$ with $\range(U_*)\subseteq\range(U)$ and
$\range(V_*)\subseteq\range(V)$.

If, in addition, $\Lambda_0\succ0$, in particular if
$\nabla^2h_i(C_0)\succeq\alpha I_{k_0}$ for every $i$ and
$C_0\in\R^{r_0\times s_0}$ and some $\alpha>0$, then for
$1\le r\le m_1$ and $1\le s\le m_2$, a lossless $(r,s)$ support
exists if and only if $r\ge r_*$ and $s\ge s_*$.
\end{corollary}

\begin{proof}[Proof of Corollary~\ref{cor:least}]
Write $G_j:=\matop(Bs_{0j})\in\R^{m_1\times m_2}$. By
\eqref{eq:feature-common-spans}, $\range(G_j)\subseteq\range(U_*)$ and
$\range(G_j^\top)\subseteq\range(V_*)$. Since $U_*$ and $V_*$ have
orthonormal columns, $U_*U_*^\top$ and $V_*V_*^\top$ are the orthogonal
projectors onto $\range(U_*)$ and $\range(V_*)$. Thus
$U_*U_*^\top G_j=G_j$ and $G_jV_*V_*^\top=G_j$, so that
\[
G_j=U_*U_*^\top G_jV_*V_*^\top=U_*(U_*^\top G_jV_*)V_*^\top.
\]
Vectorizing yields $Bs_{0j}\in\range(V_*\otimes U_*)$ for every $j$, and
hence $\range(BS_0)\subseteq\range(S_*)$. Theorem~\ref{thm:ridge}
therefore shows that $S_*$ is lossless. If $S=V\otimes U$ is an
orthonormal support with $\range(U_*)\subseteq\range(U)$ and
$\range(V_*)\subseteq\range(V)$, then $U_*=UU^\top U_*$ and
$V_*=VV^\top V_*$, so $S_*=(V\otimes U)(V^\top V_*\otimes U^\top U_*)$
and $\range(BS_0)\subseteq\range(S_*)\subseteq\range(S)$. Hence $S$ is
lossless as well. Extending $U_*$ and $V_*$ by orthonormal columns yields
lossless $(r,s)$ supports for every $r_*\le r\le m_1$ and
$s_*\le s\le m_2$.

Now assume $\Lambda_0\succ0$, and let $S=V\otimes U$ be a lossless
$(r,s)$ support. The necessity statement in Theorem~\ref{thm:ridge}
gives $\range(BS_0)\subseteq\range(S)$. Therefore each $G_j$ has a
representation $G_j=UC_jV^\top$ with $C_j\in\R^{r\times s}$, so
$\range(G_j)\subseteq\range(U)$ and $\range(G_j^\top)\subseteq\range(V)$.
Using \eqref{eq:feature-common-spans}, we obtain
\[
 \range(U_*)\subseteq\range(U),
 \qquad
 \range(V_*)\subseteq\range(V).
\]
Taking dimensions gives $r\ge r_*$ and $s\ge s_*$, proving necessity.
\end{proof}

\section{Wasserstein Budgets}
\label{sec:w-structure}

We next consider the Wasserstein budget, which constrains only the total
displacement energy rather than the full second-moment matrix. With a scalar budget, the worst case can be described by a single
response per scenario. We therefore analyze the excluded risk without the feature
structure \eqref{eq:ridge}, requiring only that each $f_i$ be twice
continuously differentiable with
\begin{equation}\label{eq:regular}
\|\nabla^2 f_i(Z)\|_{\rm op}\le L
\qquad\text{for all }i\text{ and }Z.
\end{equation}
For a prescribed radius $\rho>0$, suppose there exists $\lambda_*>L/2$ such that the maximizers
\begin{equation}\label{eq:calibrated-response}
Z_i^*:=\argmax_{Z}\{f_i(Z)-\lambda_*\|Z\|_F^2\}, \ i=1,\ldots,N, \ \text{satisfy}\ \sum_i p_i\|Z_i^*\|_F^2=\rho^2.
\end{equation}
The strict inequality $\lambda_*>L/2$ makes each penalized objective
strongly concave, so $Z_i^*$ is well defined and unique. The energy
equality calibrates these responses to the prescribed Wasserstein radius.
Lemma~\ref{lem:calibration} below shows that \eqref{eq:calibrated-response}
holds whenever the full-space dual admits an optimal multiplier above
$L/2$, and in particular for all sufficiently small radii.

\begin{theorem}[Excluded risk under Wasserstein budgets]
\label{prop:energy}
Assume \eqref{eq:regular} and \eqref{eq:calibrated-response}.
Then 
\[
R^{\rm W}(\rho)=\sum_ip_if_i(Z_i^*).
\]
Let $S=V\otimes U$ be an orthonormal support with $U\in\R^{m_1\times r}$
and $V\in\R^{m_2\times s}$, and set
$a_S:=(\sum_i p_i\|U^\top Z_i^*V\|_F^2)^{1/2}$.
Then $0\le a_S\le\rho$, and the excluded risk satisfies
\begin{equation}\label{eq:energy-cost}
(2\lambda_*-L)\rho(\rho-a_S)
\le
\Delta_S^{\rm W}
\le
(2\lambda_*+L)\rho(\rho-a_S).
\end{equation}
In particular,
\begin{equation}\label{eq:w-exact}
\Delta_S^{\rm W}=0
\quad\Longleftrightarrow\quad
UU^\top Z_i^*VV^\top=Z_i^*
\quad\text{for every }i.
\end{equation}
\end{theorem}

\begin{proof}[Proof of Theorem \ref{prop:energy}]
Let $\psi_i(Z):=f_i(Z)-\lambda_*\|Z\|_F^2$.
Since $\sum_i p_i\|Z_i^*\|_F^2=\rho^2$, the deterministic responses
$Z_i^*$ are feasible for the full-space problem. Moreover, weak duality,
the pointwise optimality of $Z_i^*$, and the energy equality give
\begin{equation}\label{eq:w-full-optimum}
\sum_i p_i f_i(Z_i^*)
\le R^{\rm W}(\rho)
\le \lambda_*\rho^2+\sum_i p_i\sup_Z\psi_i(Z)
= \lambda_*\rho^2+\sum_i p_i\psi_i(Z_i^*)
= \sum_i p_i f_i(Z_i^*).
\end{equation}
Thus the responses are full-space optimal, and $\lambda_*$ is
full-space dual-optimal.

The map $Z\mapsto UU^\top ZVV^\top$ is the orthogonal projection onto the subspace
$\{UCV^\top:C\in\R^{r\times s}\}$, and
$\|UU^\top ZVV^\top\|_F=\|U^\top ZV\|_F$. Hence
$\|Z\|_F^2=\|U^\top ZV\|_F^2+\|Z-UU^\top ZVV^\top\|_F^2$ for every $Z$.
Applying this to each $Z_i^*$ and using $\sum_ip_i\|Z_i^*\|_F^2=\rho^2$
gives
\begin{equation}\label{eq:w-omitted-energy}
\rho^2-a_S^2
=\sum_ip_i\|Z_i^*-UU^\top Z_i^*VV^\top\|_F^2\ge0,
\end{equation}
and in particular $0\le a_S\le\rho$.

We first prove the lower bound in \eqref{eq:energy-cost}. By
\eqref{eq:regular},
\[
(2\lambda_*-L)I\preceq-\nabla^2\psi_i(Z)\preceq(2\lambda_*+L)I
\qquad\text{for all }Z,
\]
and since $\nabla\psi_i(Z_i^*)=0$, Taylor's formula gives
\begin{equation}\label{eq:bregman-bound}
\Big(\lambda_*-\frac L2\Big)\|Z-Z_i^*\|_F^2
\le\psi_i(Z_i^*)-\psi_i(Z)
\le\Big(\lambda_*+\frac L2\Big)\|Z-Z_i^*\|_F^2 .
\end{equation}
Take any feasible conditional core distribution $q$, write
$Z=UCV^\top$, and let $b:=(\E_q\|C\|_F^2)^{1/2}\le\rho$.
By Cauchy--Schwarz,
\[
\E_q\langle Z,Z_i^*\rangle_F
=\E_q\langle C,U^\top Z_i^*V\rangle_F
\le(\E_q\|C\|_F^2)^{1/2}\cdot
\left(\sum_i p_i\|U^\top Z_i^*V\|_F^2\right)^{1/2}
=ba_S.
\]
Consequently, since $\|Z\|_F=\|C\|_F$ by orthonormality and
$\E_q\|Z_i^*\|_F^2=\rho^2$,
\begin{equation}\label{eq:w-distance}
\E_q\|Z-Z_i^*\|_F^2
=b^2+\rho^2-2\E_q\langle Z,Z_i^*\rangle_F
\ge b^2+\rho^2-2ba_S .
\end{equation}
By \eqref{eq:w-full-optimum} and $f_i(Z)=\psi_i(Z)+\lambda_*\|Z\|_F^2$,
\begin{equation}\label{eq:w-full-gap-identity}
R^{\rm W}(\rho)-J_S(q)
=\lambda_*(\rho^2-b^2)+\E_q[\psi_i(Z_i^*)-\psi_i(Z)].
\end{equation}
Combining the lower bound in \eqref{eq:bregman-bound} with
\eqref{eq:w-distance} yields
\[
\begin{aligned}
R^{\rm W}(\rho)-J_S(q)
&\ge\lambda_*(\rho^2-b^2)
 +\left(\lambda_*-\frac L2\right)(\rho^2+b^2-2ba_S)\\
&=(2\lambda_*-L)\rho(\rho-a_S)
 +\frac L2(\rho^2-b^2)+(2\lambda_*-L)a_S(\rho-b)\\
&\ge(2\lambda_*-L)\rho(\rho-a_S),
\end{aligned}
\]
where the last two terms on the second line are nonnegative because
$0\le b\le\rho$, $a_S\ge0$ and $2\lambda_*>L\ge0$. Taking the infimum
over feasible $q$ turns the left-hand side into $\Delta_S^{\rm W}$, which
proves the lower bound in \eqref{eq:energy-cost}.

We next prove the upper bound in \eqref{eq:energy-cost}. If $a_S=0$,
the kernel $q(dC\mid\xi_i):=\delta_0(dC)$ is feasible with $b=0$, so
\eqref{eq:w-full-gap-identity} and the upper bound in
\eqref{eq:bregman-bound} give
\[
\Delta_S^{\rm W}\le\lambda_*\rho^2+\left(\lambda_*+\frac L2\right)\sum_ip_i\|Z_i^*\|_F^2
=\left(2\lambda_*+\frac L2\right)\rho^2\le(2\lambda_*+L)\rho^2.
\]
Let $a_S>0$. We rescale the projected responses to use the full energy
budget, setting $\widehat C_i:=(\rho/a_S)\,U^\top Z_i^*V$ and
   $\widehat Z_i:=U\widehat C_iV^\top$, so that
\[
\sum_ip_i\|\widehat Z_i\|_F^2=\sum_ip_i\|\widehat C_i\|_F^2=\rho^2,
\qquad
\sum_ip_i\langle\widehat Z_i,Z_i^*\rangle_F
=\frac{\rho}{a_S}\sum_ip_i\|U^\top Z_i^*V\|_F^2=\rho a_S,
\]
and hence
\[
\sum_ip_i\|\widehat Z_i-Z_i^*\|_F^2=2\rho^2-2\rho a_S=2\rho(\rho-a_S).
\]
The conditional distribution
$\widehat q(dC\mid\xi_i):=\delta_{\widehat C_i}(dC)$ is feasible because
$\sum_ip_i\|\widehat C_i\|_F^2=\rho^2$, so
$J_S(\widehat q)\le R_S^{\rm W}(\rho)$. Moreover $b=\rho$ for
$\widehat q$, so the first term in \eqref{eq:w-full-gap-identity}
vanishes. Therefore, by the upper bound in \eqref{eq:bregman-bound},
\[
\begin{aligned}
\Delta_S^{\rm W}
\le R^{\rm W}(\rho)-J_S(\widehat q)
=\sum_ip_i[\psi_i(Z_i^*)-\psi_i(\widehat Z_i)]
&\le\Big(\lambda_*+\frac L2\Big)\sum_ip_i\|\widehat Z_i-Z_i^*\|_F^2\\
&=(2\lambda_*+L)\rho(\rho-a_S),
\end{aligned}
\]
which proves the upper bound in \eqref{eq:energy-cost}.

Finally, since $2\lambda_*-L>0$ and $\rho>0$, \eqref{eq:energy-cost}
gives $\Delta_S^{\rm W}=0$ if and only if $a_S=\rho$. By
\eqref{eq:w-omitted-energy} and $p_i>0$, this is equivalent to
$UU^\top Z_i^*VV^\top=Z_i^*$ for every $i$, proving \eqref{eq:w-exact}.
\end{proof}

The proof shows that $\lambda_*$ in \eqref{eq:calibrated-response} is a
dual-optimal multiplier for the full-space problem. Conversely, the following lemma shows that
if $\lambda_*>L/2$ is full-space dual-optimal, its pointwise
maximizing responses satisfy the energy equality in
\eqref{eq:calibrated-response}. It also gives an a priori
sufficient condition on $\rho$ ensuring $\lambda_*>L/2$.

\begin{lemma}[Calibration and a priori radius condition]\label{lem:calibration}
Assume \eqref{eq:regular}, and let $\lambda_*\ge0$ be a dual-optimal
multiplier of $R^{\rm W}(\rho)$.
If $\lambda_*>L/2$, then the maximizers
$Z_i^*$ of $f_i(Z)-\lambda_*\|Z\|_F^2$ satisfy
$\sum_ip_i\|Z_i^*\|_F^2=\rho^2$. Moreover, if $L>0$ and
\begin{equation}\label{eq:calibration-apriori}
0<\rho<\frac{1}{2L}\left(\sum_ip_i\|\nabla f_i(0)\|_F^2\right)^{1/2},
\end{equation}
then every dual-optimal multiplier satisfies $\lambda_*>L/2$.
\end{lemma}

\begin{proof}[Proof of Lemma \ref{lem:calibration}]
Write $\mathcal D^{\rm W}$ for the dual function \eqref{eq:dual-value-w} of the
full-space problem, i.e., with $S=I_d$. Let $\lambda_*>L/2$ be dual
optimal and $\psi_i(Z):=f_i(Z)-\lambda_*\|Z\|_F^2$. By \eqref{eq:regular},
$\psi_i$ is $(2\lambda_*-L)$-strongly concave with unique maximizer
$Z_i^*$, so
\[
\psi_i(Z_i^*)-\psi_i(Z)\ge\left(\lambda_*-\frac L2\right)\|Z-Z_i^*\|_F^2
\qquad\text{for all }Z.
\]
By \eqref{eq:strong-duality}, $R^{\rm W}(\rho)=\mathcal D^{\rm W}(\lambda_*)$, so
there are feasible conditional distributions $q_m$ whose expected losses
converge to $\mathcal D^{\rm W}(\lambda_*)$, and hence $s(q_m)+e(q_m)\to0$ by
\eqref{eq:gap-identity} with multiplier $\lambda_*$. The preceding
inequality and \eqref{eq:response-error} give
\[
e(q_m)\ge(\lambda_*-L/2)\,\E_{q_m}\|Z-Z_i^*\|_F^2,
\]
so
$\E_{q_m}\|Z-Z_i^*\|_F^2\to0$, while \eqref{eq:budget-slack} and
$\lambda_*>0$ give $\E_{q_m}\|Z\|_F^2\to\rho^2$. By Minkowski's inequality,
\[
\left|\left(\E_{q_m}\|Z\|_F^2\right)^{1/2}
-\left(\sum_ip_i\|Z_i^*\|_F^2\right)^{1/2}\right|
\le\left(\E_{q_m}\|Z-Z_i^*\|_F^2\right)^{1/2}\to0,
\]
so $\sum_ip_i\|Z_i^*\|_F^2=\lim_m\E_{q_m}\|Z\|_F^2=\rho^2$, which proves
\eqref{eq:calibrated-response}.

Now assume $L>0$ and \eqref{eq:calibration-apriori}. For $\mu>L/2$, let
$Z_i(\mu)$ be the unique maximizer of $f_i(Z)-\mu\|Z\|_F^2$. Then
$\nabla f_i(Z_i(\mu))=2\mu Z_i(\mu)$, and \eqref{eq:regular} gives
$\|\nabla f_i(0)\|_F\le\|\nabla f_i(Z_i(\mu))\|_F+L\|Z_i(\mu)\|_F
=(2\mu+L)\|Z_i(\mu)\|_F$. Hence
\[
\sum_ip_i\|Z_i(\mu)\|_F^2
\ge\frac{\sum_ip_i\|\nabla f_i(0)\|_F^2}{(2\mu+L)^2}.
\]
As $\mu\downarrow L/2$ the right-hand side tends to
$(2L)^{-2}\sum_ip_i\|\nabla f_i(0)\|_F^2$, which exceeds $\rho^2$ by
\eqref{eq:calibration-apriori}. Fix such a $\mu>L/2$ with
$\sum_ip_i\|Z_i(\mu)\|_F^2>\rho^2$. For $0\le\lambda<\mu$,
evaluating the supremum in $\mathcal D^{\rm W}(\lambda)$ at $Z_i(\mu)$ gives
\[
\mathcal D^{\rm W}(\lambda)
\ge\lambda\rho^2+\sum_ip_i\bigl[f_i(Z_i(\mu))-\lambda\|Z_i(\mu)\|_F^2\bigr]
=\mathcal D^{\rm W}(\mu)+(\mu-\lambda)\left(\sum_ip_i\|Z_i(\mu)\|_F^2-\rho^2\right)
>\mathcal D^{\rm W}(\mu).
\]
Thus no $\lambda<\mu$ minimizes $\mathcal D^{\rm W}$, and every dual-optimal
multiplier satisfies $\lambda_*\ge\mu>L/2$.
\end{proof}

Condition \eqref{eq:w-exact} likewise determines how many common column
and row directions a lossless support needs and how to choose them.

\begin{corollary}[Minimal lossless dimensions]
\label{cor:w-least}
Under the assumptions of Theorem~\ref{prop:energy}, let
$U_*\in\R^{m_1\times r_*}$ and $V_*\in\R^{m_2\times s_*}$ have
orthonormal columns with
\begin{equation}\label{eq:w-response-spans}
 \range(U_*)=\sum_{i=1}^{N}\range(Z_i^*),
 \qquad
 \range(V_*)=\sum_{i=1}^{N}\range(Z_i^{*\top}),
\end{equation}
where $r_*$ and $s_*$ are the dimensions of these two subspaces. Then an
orthonormal support $S=V\otimes U$ is lossless if and only if
\begin{equation}\label{eq:w-lossless-factor-inclusions}
\range(U_*)\subseteq\range(U),
\qquad
\range(V_*)\subseteq\range(V).
\end{equation}
In particular, $S_*:=V_*\otimes U_*$ is lossless, and for
$1\le r\le m_1$ and $1\le s\le m_2$, a lossless $(r,s)$ support exists
if and only if $r\ge r_*$ and $s\ge s_*$.
\end{corollary}

\begin{proof}[Proof of Corollary~\ref{cor:w-least}]
Let $U$ and $V$ have orthonormal columns, and set $S=V\otimes U$.
Since $UU^\top$ and $VV^\top$ are orthogonal projectors, for every
$Z\in\R^{m_1\times m_2}$,
\[
\begin{aligned}
 UU^\top ZVV^\top=Z
 &\quad\Longleftrightarrow\quad
 UU^\top Z=Z,\quad ZVV^\top=Z\\
 &\quad\Longleftrightarrow\quad
 \range(Z)\subseteq\range(U),\quad
 \range(Z^\top)\subseteq\range(V).
\end{aligned}
\]
Applying this to each $Z_i^*$ and using
\eqref{eq:w-response-spans}, the losslessness condition
\eqref{eq:w-exact} is equivalent to
\eqref{eq:w-lossless-factor-inclusions}.

Since $\sum_i p_i\|Z_i^*\|_F^2=\rho^2>0$, at least one $Z_i^*$
is nonzero, so $r_*,s_*\ge1$.
The support $S_*=V_*\otimes U_*$ satisfies
\eqref{eq:w-lossless-factor-inclusions} and is therefore lossless.
For any admissible $r\ge r_*$ and $s\ge s_*$, extending the
orthonormal bases $U_*$ and $V_*$ yields a lossless $(r,s)$ support.
Conversely, every lossless $(r,s)$ support satisfies
\eqref{eq:w-lossless-factor-inclusions}, which implies
$r\ge r_*$ and $s\ge s_*$.
\end{proof}

\begin{remark}[Truncated supports]
\label{rem:w-truncation}
The subspaces in \eqref{eq:w-response-spans} can be the entire ambient
spaces, so a low-dimensional lossless support need not exist. When the
prescribed dimensions $(r,s)$ admit no lossless support, an
energy-ordered truncation of the bases in \eqref{eq:w-response-spans}
gives a support whose excluded risk is within an explicit constant
factor of the best achievable.
Specifically, choose the columns of $U_*$ and $V_*$ as orthonormal
eigenvectors corresponding to the positive eigenvalues of
\[
 \sum_i p_i Z_i^*Z_i^{*\top}
 \quad\text{and}\quad
 \sum_i p_i Z_i^{*\top}Z_i^*,
\]
respectively, ordered by decreasing eigenvalue.
Let $\widehat U$ and $\widehat V$ consist of the first $r$ and $s$
of these columns, respectively, extending $U_*$ or $V_*$ by
orthonormal columns if $r>r_*$ or $s>s_*$, and set
$\widehat S:=\widehat V\otimes\widehat U$.
The matrices $\widehat U$ and $\widehat V$ minimize the column and row
projection errors of the responses $Z_i^*$ separately. For any
$(r,s)$ support $S=V\otimes U$, since
$\|U^\top Z_i^*V\|_F\le\min\{\|U^\top Z_i^*\|_F,\|Z_i^*V\|_F\}$, each of
these minimal errors is at most the total omitted response energy
$\rho^2-a_S^2$. Together with the orthogonal decomposition
$Z-UU^\top ZVV^\top=(I-UU^\top)Z+UU^\top Z(I-VV^\top)$,
this gives
\[
 \rho^2-a_{\widehat S}^2
 \le
 \sum_i p_i\|(I-\widehat U\widehat U^\top)Z_i^*\|_F^2
 +\sum_i p_i\|Z_i^*(I-\widehat V\widehat V^\top)\|_F^2
 \le 2(\rho^2-a_S^2).
\]
Combining this with \eqref{eq:energy-cost} and
$\rho(\rho-a_S)\le\rho^2-a_S^2\le2\rho(\rho-a_S)$ yields
\[
 \Delta_{\widehat S}^{\rm W}
 \le
 \frac{4(2\lambda_*+L)}{2\lambda_*-L}
 \inf_S\Delta_S^{\rm W},
\]
where the infimum is over all $(r,s)$ supports.
\end{remark}

\section{Worst-Case Generation}\label{sec:recovery}

 This section
computes a feasible conditional distribution whose risk is close to $R_S$.
Fix the decision $\theta$ and the shared support $S=V\otimes U$. 
We optimize over deterministic response maps $T$, with
$T(\xi_i)=T_i\in\R^{r\times s}$ assigning a core perturbation to each nominal
observation. Since $J_S$ need not be concave, we introduce a fixed
curvature correction whose maximization is concave and upper bounds $R_S$.
This is equivalent to restricting the
original dual multipliers from below by the curvature correction
(Section~\ref{sec:map-formulation}).
Projected gradient ascent solves the corrected problem while keeping every
iterate within the original budget, and at every iterate it reports an
interval that contains $R_S$ (Section~\ref{sec:global-projected-generation}).
The interval closes exactly when the correction vanishes at the limit, and a
strict multiplier margin yields linear convergence
(Section~\ref{sec:exact-generation}). Adding the excluded risk $\Delta_S$ in
\eqref{eq:total-error}, which Theorems~\ref{thm:ridge} and~\ref{prop:energy}
control, to the width of this interval bounds the gap to the unrestricted
risk $R$.

\subsection{Deterministic responses and curvature correction}
\label{sec:map-formulation}

A deterministic response
assigns a core matrix $T_i\in\R^{r\times s}$ to each atom. For
$T=(T_i)_{i=1}^N$, define
\[
 q_T(dC\mid\xi_i)=\delta_{T_i}(dC),\qquad
 J_S(T):=J_S(q_T),\qquad M(T):=M(q_T),
\]
and $\|T\|_{L^2}^2:=\sum_i p_i\|T_i\|_F^2$.
The Lagrangians evaluated at this distribution are
\begin{align*}
 \mathcal L(\theta,T;\Lambda)
 &=J_S(T)+\tr[\Lambda(B_S-M(T))],
 \qquad\Lambda\succeq0,\\
 \mathcal L(\theta,T;\lambda)
 &=J_S(T)+\lambda(\rho^2-\|T\|_{L^2}^2),
 \qquad\lambda\ge0.
\end{align*}
The suprema in \eqref{eq:dual-value} and \eqref{eq:dual-value-w} are taken
separately at each atom, so maximizing these Lagrangians over deterministic
responses recovers the dual functions:
\[
 \mathcal D_S^{\rm M}(\Lambda)=\sup_T\mathcal L(\theta,T;\Lambda),\qquad
 \mathcal D_S^{\rm W}(\lambda)=\sup_T\mathcal L(\theta,T;\lambda).
\]
Consequently, if a feasible $T$ and a dual-optimal multiplier $\Lambda$
satisfy $\mathcal L(\theta,T;\Lambda)=\mathcal D_S^{\rm M}(\Lambda)$ and
$\tr[\Lambda(B_S-M(T))]=0$, then $J_S(T)=\mathcal D_S^{\rm M}(\Lambda)=R_S$ by
Proposition~\ref{prop:duality}, and likewise for the Wasserstein budget.
Such a response need not exist, since the complementarity condition may
require randomizing among several maximizers.
We therefore construct responses from a concave upper problem and certify
their gap to $R_S$.

Throughout this section, each $f_i$ is twice continuously differentiable and
satisfies the curvature bound \eqref{eq:regular} of
Section~\ref{sec:w-structure}, which also holds under \eqref{eq:ridge}.
By Taylor's formula, this bound implies Assumption~\ref{ass:growth} with
$a(\xi)=\|\nabla f_\xi(0)\|_F$.
Write $c_i:=\vecop(T_i)$ and let
$f_i^S(c):=f_i(U\matop(c)V^\top)$
denote the loss restricted to the chosen support. We use the matrix-shaped
core gradient
\[
g_i(T_i):=\matop(\nabla f_i^S(c_i))=U^\top\nabla f_i(UT_iV^\top)V.
\]
By the chain rule,
$\nabla^2 f_i^S(c)=S^\top\nabla^2 f_i(U\matop(c)V^\top)S$ in vectorized
coordinates, and since $S^\top S=I_k$, \eqref{eq:regular} gives
\begin{equation}\label{eq:restricted-curvature}
 -LI_k\preceq\nabla^2 f_i^S(c)\preceq LI_k
 \qquad\text{for all } i \text{ and } c.
\end{equation}
Fix a positive semidefinite matrix $\Omega\in\R^{k\times k}$ such that
\begin{equation}\label{eq:curvature-correction}
 \nabla^2 f_i^S(c)\preceq2\Omega
 \qquad\text{for all } i \text{ and } c.
\end{equation}
By \eqref{eq:restricted-curvature}, $\Omega=\frac L2I_k$ is always
admissible, and a smaller $\Omega$ in the Loewner order yields a smaller
correction below and a weaker multiplier restriction in
Proposition~\ref{prop:exact-shift-recovery}, which matters when the
curvature of the loss is anisotropic. 
The feasible response sets are
\[
 \mathcal K_{\rm M}:=\{T:M(T)\preceq B_S\},\qquad
 \mathcal K_{\rm W}:=\{T:\|T\|_{L^2}\le\rho\}.
\]
Both sets are compact and convex.
Define the corrected objectives
\[
 \Phi_{\rm M}(T)
 :=J_S(T)+\tr[\Omega(B_S-M(T))],\qquad
 \Phi_{\rm W}(T)
 :=J_S(T)+\frac L2(\rho^2-\|T\|_{L^2}^2),
\]
that is, $\Phi_{\rm M}(T)=\mathcal L(\theta,T;\Omega)$ and
$\Phi_{\rm W}(T)=\mathcal L(\theta,T;L/2)$. 
Equivalently,
\[
 \Phi_{\rm M}(T)=\tr(\Omega B_S)
 +\sum_ip_i\big[f_i^S(c_i)-c_i^\top\Omega c_i\big],
 \qquad
 \Phi_{\rm W}(T)=\frac L2\rho^2
 +\sum_ip_i\Big[f_i^S(c_i)-\frac L2\|c_i\|^2\Big],
\]
and the summands have Hessians $\nabla^2f_i^S(c)-2\Omega\preceq0$ and
$\nabla^2f_i^S(c)-LI_k\preceq0$ by \eqref{eq:curvature-correction} and
\eqref{eq:restricted-curvature}. Hence both objectives are concave.
For the chosen budget, write $\Phi$ and $\mathcal K$ for the corresponding
objective and feasible set, and define
$\Phi^\star:=\max_{T\in\mathcal K}\Phi(T)$.

\begin{proposition}[Concave surrogate and restricted duality]
\label{prop:exact-shift-recovery}
Under \eqref{eq:regular}, for every $T\in\mathcal K$,
\[
 J_S(T)\le R_S\le\Phi^\star
 \quad\text{and}\quad
 J_S(T)\le\Phi(T).
\]
Moreover,
\begin{align}
 \max_{T\in\mathcal K_{\rm M}}\Phi_{\rm M}(T)
  &=\min_{\Lambda\succeq\Omega}\mathcal D_S^{\rm M}(\Lambda),
 \label{eq:restricted-dual-matrix}\\
 \max_{T\in\mathcal K_{\rm W}}\Phi_{\rm W}(T)
 &=\min_{\lambda\ge L/2}\mathcal D_S^{\rm W}(\lambda).
 \label{eq:restricted-dual-scalar}
\end{align}
Every maximizer $\bar T$ can be paired with an optimal multiplier
of the corresponding restricted dual satisfying the following KKT
conditions. In the moment case,
\begin{equation}
\label{eq:corrected-matrix-kkt}
    \vecop(g_i(\bar T_i))=2\bar\Lambda\vecop(\bar T_i),
 \ i=1,\ldots,N,\quad
  \bar\Lambda\succeq\Omega,\quad
 \tr[(\bar\Lambda-\Omega)(B_S-M(\bar T))]=0.
\end{equation}
In the Wasserstein case,
\begin{equation}
\label{eq:corrected-scalar-kkt}
    g_i(\bar T_i)=2\bar\lambda\bar T_i,
 \ i=1,\ldots,N,\quad \bar\lambda\ge L/2,\quad
 (\bar\lambda-L/2)(\rho^2-\|\bar T\|_{L^2}^2)=0.
\end{equation}

\end{proposition}

\begin{proof}[Proof of Proposition \ref{prop:exact-shift-recovery}]

Let $q$ be feasible and define its conditional mean response by
$
T_i:=\E_q[C\mid \xi_i]$, $i=1,\ldots,N$.
By conditional Jensen's inequality,
\[
M(T)\preceq M(q),
\qquad
\|T\|_{L^2}^2\le \E_q\|C\|_F^2,
\]
so $T\in\mathcal K$ under either budget. Since
$c\mapsto f_i^S(c)-c^\top\Omega c$ and $c\mapsto f_i^S(c)-L\|c\|^2/2$ are
concave,
\[
J_S(q)
\le
J_S(T)+\tr[\Omega(M(q)-M(T))]
\qquad\text{and}\qquad
J_S(q)
\le
J_S(T)
+\frac L2\bigl(\E_q\|C\|_F^2-\|T\|_{L^2}^2\bigr).
\]
Under the moment budget, $\Omega\succeq0$ and $M(q)\preceq B_S$ give
$\tr[\Omega M(q)]\le\tr[\Omega B_S]$, whereas under the
Wasserstein budget, $\E_q\|C\|_F^2\le\rho^2$. Hence, in either case,
$
J_S(q)\le \Phi(T)$.
Taking the supremum over feasible $q$ gives $R_S\le\Phi^\star$.
Every $T\in\mathcal K$ defines a feasible deterministic response, so
$J_S(T)\le R_S$. Moreover, since $\Omega\succeq0$ and $L\ge0$, the correction is nonnegative
on $\mathcal K$, and hence $J_S(T)\le\Phi(T)$. Also, since $\Phi$ is continuous and
$\mathcal K$ is compact, $\Phi^\star$ is attained.

For the moment problem, introduce a multiplier $\Theta\succeq0$ for the
constraint $M(T)\preceq B_S$. Its Lagrangian is
\[
\Phi_{\rm M}(T)+\tr[\Theta(B_S-M(T))]
=
\mathcal L(\theta,T;\Omega+\Theta).
\]
Hence, the dual
function is
$
\mathcal D_S^{\rm M}(\Omega+\Theta)$.
The problem is a finite-dimensional concave maximization with a
matrix convex constraint. Moreover,
$T=0$ is strictly feasible since $M(0)=0\prec B_S$. Slater's condition
therefore yields strong duality and dual attainment. With
$
\Lambda=\Omega+\Theta$,
the dual constraint $\Theta\succeq0$ becomes
$\Lambda\succeq\Omega$, which proves
\eqref{eq:restricted-dual-matrix}.
The KKT stationarity and complementary
slackness conditions give \eqref{eq:corrected-matrix-kkt}.
For the Wasserstein problem, introducing $\tau\ge0$ for
$\|T\|_{L^2}^2\le\rho^2$ gives the Lagrangian
$
\mathcal L(\theta,T;L/2+\tau)$.
Again $T=0$ is strictly feasible, and Slater's condition yields strong
duality and dual attainment. Setting $\lambda=L/2+\tau$ gives
\eqref{eq:restricted-dual-scalar}, while the corresponding KKT conditions
give \eqref{eq:corrected-scalar-kkt}.
\end{proof}

The bound $R_S\le\Phi^\star$ concerns the original problem over
conditional distributions, not only deterministic responses: by the
concavity of the corrected losses, every feasible $q$ satisfies
$J_S(q)\le\Phi(T)$ for its conditional mean response $T$, which remains
feasible.
The correction also has a dual interpretation: it restricts the original
dual multiplier to $\Lambda\succeq\Omega$ under the moment budget and to
$\lambda\ge L/2$ under the Wasserstein budget.
Since both the original and restricted dual minima are attained,
\eqref{eq:strong-duality}, \eqref{eq:restricted-dual-matrix}, and
\eqref{eq:restricted-dual-scalar} imply that
$\Phi^\star=R_S$ if and only if the original dual admits an optimal
multiplier satisfying the corresponding restriction.
Finally, every maximizer $\bar T\in\argmax_{T\in\mathcal K}\Phi(T)$
satisfies $\Phi^\star-J_S(\bar T)=\tr[\Omega(B_S-M(\bar T))]$ under the
moment budget and
$\Phi^\star-J_S(\bar T)=\frac L2(\rho^2-\|\bar T\|_{L^2}^2)$ under the
Wasserstein budget. Since $J_S(\bar T)\le R_S\le\Phi^\star$, we have
$J_S(\bar T)=R_S=\Phi^\star$ if and only if this correction vanishes.

\subsection{Projected gradient method and risk certificates}
\label{sec:global-projected-generation}

For response collections $T=(T_i)_{i=1}^N$ and $T'=(T_i')_{i=1}^N$,
we use the inner product
\begin{equation}\label{eq:response-inner-product}
 \langle T,T'\rangle=
 \begin{cases}
 \displaystyle\sum_i p_i\vecop(T_i)^\top B_S^{-1}\vecop(T_i'),
  &\text{(Moment)},\\[2pt]
 \displaystyle\sum_i p_i\langle T_i,T_i'\rangle_F,
  &\text{(Wasserstein)}.
 \end{cases}
\end{equation}
Its induced norm satisfies
$\|T\|^2=\tr[B_S^{-1}M(T)]$ in the matrix-moment case and
$\|T\|^2=\|T\|_{L^2}^2$ in the Wasserstein case.
In particular, every feasible response satisfies $\|T\|\le \rho_{\mathcal K}$,
where $\rho_{\mathcal K}=\sqrt{k}$ or $\rho_{\mathcal K}=\rho$, respectively. 
 The gradients are
\begin{equation}\label{eq:response-gradient}
 (\nabla\Phi(T))_i=
 \begin{cases}
  \matop\!\left(B_S[\vecop(g_i(T_i))-2\Omega\vecop(T_i)]\right),
  &\text{(Moment)},\\
 g_i(T_i)-LT_i, &\text{(Wasserstein)}.
 \end{cases}
\end{equation}
By \eqref{eq:restricted-curvature} and \eqref{eq:curvature-correction},
$\Phi$ is concave and has
$L_\Phi$-Lipschitz gradient with respect to the inner product
\eqref{eq:response-inner-product}, where
\[
 L_\Phi=
 \begin{cases}
 \lambda_{\max}\bigl(B_S^{1/2}(2\Omega+LI_k)B_S^{1/2}\bigr), &\text{(Moment)},\\
 2L, &\text{(Wasserstein)}.
 \end{cases}
\]
Here, in the moment case, for any response direction $T'$ with
$c_i'=\vecop(T_i')$,
\[
 \begin{aligned}
 0&\le-\langle T',\nabla^2\Phi_{\rm M}(T)[T']\rangle
 =\sum_i p_i c_i'^\top[2\Omega-\nabla^2 f_i^S(c_i)]c_i'
 \le\sum_i p_i c_i'^\top(2\Omega+LI_k)c_i'
 \le L_{\Phi}\|T'\|^2.
 \end{aligned}
\]
The Wasserstein bound follows directly from
\eqref{eq:restricted-curvature} under the weighted Frobenius geometry.
Then projected gradient ascent takes the form
\begin{equation}\label{eq:projected-generation}
 \widetilde T_{i,n+1}=T_{i,n}+\eta(\nabla\Phi(T_n))_i,
 \qquad T_{n+1}=\Pi_{\mathcal K}(\widetilde T_{n+1}),
\end{equation}
where $\Pi_{\mathcal K}(T)=\argmin_{T'\in\mathcal K}\|T'-T\|^2$.

\paragraph{Closed-form projection.} For the moment budget, let
\[
B_S^{-1/2}M(T)B_S^{-1/2}
=
Q\diag(d_1,\ldots,d_k)Q^\top
\]
be an eigendecomposition, where $Q\in\mathbb R^{k\times k}$ is
orthogonal and $d_j\ge0$. 
Then
\begin{equation}\label{eq:matrix-budget-projection}
(\Pi_{\mathcal K_{\rm M}}(T))_i
=
\matop(
B_S^{1/2}Q\diag(a_1,\ldots,a_k)Q^\top
B_S^{-1/2}\vecop(T_i)),
\end{equation}
where 
$a_j:=1/\max\{1,\sqrt{d_j}\}$, $j=1,\ldots,k$.
Actually, the map sending $T$ to the $k\times N$ matrix whose $i$th column is
$\sqrt{p_i}B_S^{-1/2}\vecop(T_i)$ is an isometry from the norm induced by
\eqref{eq:response-inner-product} to the Frobenius norm. Since this
matrix multiplied by its transpose equals $B_S^{-1/2}M(T)B_S^{-1/2}$, the
map sends $\mathcal K_{\rm M}$ onto the unit operator-norm ball. The
Frobenius projection onto this ball clips singular values at one, which
multiplies each column by $Q\diag(a_1,\ldots,a_k)Q^\top$ and gives
\eqref{eq:matrix-budget-projection}.
For the Wasserstein budget, the projection is radial:
\begin{equation}\label{eq:scalar-budget-projection}
 (\Pi_{\mathcal K_{\rm W}}(T))_i
 =\frac{\rho}{\max\{\rho,\|T\|_{L^2}\}}T_i.
\end{equation}

\paragraph{Risk certificates.}\label{sec:feasible-certificate}
For a feasible response $T\in\mathcal K$, define the first-order
optimality gap
\[
\mathcal G(T)
:=
\max_{T'\in\mathcal K}
\langle\nabla\Phi(T),T'-T\rangle.
\]
It measures the largest feasible increase of the affine first-order model
of $\Phi$ at $T$. 
The gap has the closed forms
\begin{equation}\label{eq:projected-gap-closed-forms}
 \mathcal G(T)=
 \begin{cases}
 \tr\left[
 (B_S^{-1/2}M(\nabla\Phi(T))B_S^{-1/2})^{1/2}\right]
 -\langle\nabla\Phi(T),T\rangle, &\text{(Moment)},\\
 \rho\|\nabla\Phi(T)\|-\langle\nabla\Phi(T),T\rangle, &\text{(Wasserstein)}.
 \end{cases}
\end{equation}
In the Wasserstein case, the maximum of a
linear function over the ball of radius $\rho$ is $\rho$ times the norm of
its gradient. In the moment case, it is the maximum of a linear function over the unit
operator-norm ball, which equals the nuclear norm, i.e., the sum of the
singular values, of the corresponding matrix.
We use it to form the computable upper certificate
\[
\bar R(T):=\Phi(T)+\mathcal G(T).
\]

\begin{proposition}[Computable risk certificate]
\label{prop:first-order-certificate}
Under \eqref{eq:regular}, every $T\in\mathcal K$ satisfies
\[
0\le
\Phi^\star-\Phi(T)
\le
\mathcal G(T).
\]
In particular, $\mathcal G(T)=0$ if and only if
$T\in\argmax_{T'\in\mathcal K}\Phi(T')$. Moreover,
\[
 J_S(T)\le R_S\le\Phi^\star\le\bar R(T).
\]
\end{proposition}

\begin{proof}[Proof of Proposition \ref{prop:first-order-certificate}]
By concavity, for every $T'\in\mathcal K$,
\[
\Phi(T')
\le
\Phi(T)
+\langle\nabla\Phi(T),T'-T\rangle.
\]
Maximizing over $T'$ gives
$\Phi^\star-\Phi(T)\le\mathcal G(T)$, while
$\mathcal G(T)\ge0$ follows by choosing $T'=T$.
The characterization $\mathcal G(T)=0$ follows from the first-order
optimality condition for concave maximization over the convex set
$\mathcal K$. Proposition~\ref{prop:exact-shift-recovery} gives
$J_S(T)\le R_S\le\Phi^\star$, and the preceding bound gives
$\Phi^\star\le\Phi(T)+\mathcal G(T)=\bar R(T)$.
\end{proof}

The preceding construction leads to the following projected gradient scheme with a computable stopping certificate.

\begin{algorithm}[H]
\caption{Projected gradient method with risk certification}
\label{alg:regular-generation}
\begin{algorithmic}[1]
\Statex \textbf{Input:} Fixed $\theta,U,V$, tolerance
$\varepsilon_{\rm opt}>0$, stepsize $\eta>0$, and initial response $T_0=0$.
\For{$n=0,1,\ldots$}
 \State Evaluate $J_S(T_n)$, $\Phi(T_n)$, and $\nabla\Phi(T_n)$.
 \State Compute $\mathcal G(T_n)$ from
        \eqref{eq:projected-gap-closed-forms} and set
        $\bar R(T_n)=\Phi(T_n)+\mathcal G(T_n)$.
 \If{$\mathcal G(T_n)\le\varepsilon_{\rm opt}$}
  \State \Return $T_n$ and the interval
  $[J_S(T_n),\bar R(T_n)]$.
 \EndIf
 \State Update $T_{n+1}$ by \eqref{eq:projected-generation}, using the projection
 \eqref{eq:matrix-budget-projection} or
 \eqref{eq:scalar-budget-projection}.
\EndFor
\end{algorithmic}
\end{algorithm}

\begin{theorem}[Convergence of projected gradient method]
\label{thm:generator-first-order}
Suppose \eqref{eq:regular} holds, and let $\{T_n\}_{n\ge0}$ be generated
by \eqref{eq:projected-generation} with $T_0\in\mathcal K$ and a constant
stepsize $\eta>0$ satisfying $\eta L_\Phi\le1$.
Then $T_n\in\mathcal K$ for all $n\ge0$, and there exists
$\bar T\in\argmax_{T\in\mathcal K}\Phi(T)$ such that
\[
T_n\to\bar T,
\qquad
\Phi(T_n)\uparrow\Phi^\star.
\]
Moreover, for every $T^\star\in\argmax_{T\in\mathcal K}\Phi(T)$ and
every $n\ge1$,
\begin{align}
\label{eq:projected-sublinear}
0\le\Phi^\star-\Phi(T_n)\le
\frac{\|T_0-T^\star\|^2}{2\eta n},\qquad
0\le\mathcal G(T_n)\le
\frac{4\rho_{\mathcal K}\|T_0-T^\star\|}{\eta n}.
\end{align}
Consequently, $\bar R(T_n)\to\Phi^\star$, and
\begin{equation}\label{eq:limiting-certificate-width}
\lim_{n\to\infty}
\bigl[\bar R(T_n)-J_S(T_n)\bigr]
=
\Phi^\star-J_S(\bar T)
=
\begin{cases}
\tr[\Omega(B_S-M(\bar T))],
&\text{(Moment)},\\[2mm]
\dfrac L2(\rho^2-\|\bar T\|_{L^2}^2),
&\text{(Wasserstein)}.
\end{cases}
\end{equation}
\end{theorem}

\begin{proof}[Proof of Theorem \ref{thm:generator-first-order}]
Recall that the norms and inner products are those of
\eqref{eq:response-inner-product}.
Feasibility follows from the projection step.
The characterization of the convex projection gives
\[
\eta\langle\nabla\Phi(T_n),T'-T_{n+1}\rangle
\le
\langle T_{n+1}-T_n,T'-T_{n+1}\rangle
\quad\text{for any}\ T'\in\mathcal K.
\]
Combining this inequality with concavity and $L_\Phi$-smoothness of
$\Phi$ yields
\begin{equation}\label{eq:projected-three-point}
\begin{aligned}
2\eta[\Phi(T')-\Phi(T_{n+1})]
&\le
\|T_n-T'\|^2-\|T_{n+1}-T'\|^2-(1-\eta L_\Phi)\|T_{n+1}-T_n\|^2.
\end{aligned}
\end{equation}
Taking $T'=T_n$, we obtain
\begin{equation}\label{eq:projected-ascent}
\Phi(T_{n+1})-\Phi(T_n)
\ge
\left(\frac1\eta-\frac{L_\Phi}{2}\right)
\|T_{n+1}-T_n\|^2
\ge
\frac{\|T_{n+1}-T_n\|^2}{2\eta}.
\end{equation}
Thus $\{\Phi(T_n)\}$ is nondecreasing.

Fix $T^\star\in\argmax_{T\in\mathcal K}\Phi(T)$.
Since $\eta L_\Phi\le1$, taking $T'=T^\star$ in
\eqref{eq:projected-three-point} gives
\[
2\eta[\Phi^\star-\Phi(T_{n+1})]
\le
\|T_n-T^\star\|^2-\|T_{n+1}-T^\star\|^2.
\]
Summing over the first $n$ iterations yields
\begin{equation}\label{eq:sum_obj}
2\eta\sum_{j=1}^{n}[\Phi^\star-\Phi(T_j)]
\le
\|T_0-T^\star\|^2-\|T_n-T^\star\|^2
\le
\|T_0-T^\star\|^2.
\end{equation}
Monotonicity of $\{\Phi(T_n)\}$ now proves the first part of
\eqref{eq:projected-sublinear}.

Moreover, since $\Phi(T_{n+1})\le\Phi^\star$ and
$\eta L_\Phi\le1$, the same inequality with $T'=T^\star$ gives
\[
\|T_{n+1}-T^\star\|^2
\le
\|T_n-T^\star\|^2.
\]
Hence $\{T_n\}$ is Fej\'er monotone with respect to the solution set
$\argmax_{T\in\mathcal K}\Phi(T)$.
Since $\mathcal K$ is compact, the sequence has an accumulation point
$\bar T\in\mathcal K$. By \eqref{eq:projected-sublinear} and continuity,
$\Phi(\bar T)=\Phi^\star$.
The sequence $\{\|T_n-\bar T\|\}$ is therefore nonincreasing and has
a subsequence converging to zero. Hence $T_n\to\bar T$.

We next prove the second part of
\eqref{eq:projected-sublinear}.
The Hessian bound
$-L_\Phi I\preceq\nabla^2\Phi\preceq0$, in the response inner product,
implies
\begin{equation}\label{eq:forward-nonexpansion}
\|(T+\eta\nabla\Phi(T))-(T'+\eta\nabla\Phi(T'))\|=\left\|\left[
I+\eta\int_0^1
\nabla^2\Phi\bigl(T'+t(T-T')\bigr)\,dt
\right](T-T')\right\|
\le
\|T-T'\|,
\end{equation}
where the inequality holds because the operator in brackets is
self-adjoint with spectrum contained in
$[1-\eta L_\Phi,1]\subseteq[0,1]$.
Since the metric projection is also nonexpansive, the projected gradient
map is nonexpansive. Consequently,
\[
\|T_{j+1}-T_j\|
\le
\|T_j-T_{j-1}\|,
\qquad j\ge1.
\]

For $n\ge2$ and $1\le j<n$, summing
\eqref{eq:projected-ascent} and using the preceding monotonicity gives
\[
\begin{aligned}
\Phi^\star-\Phi(T_j)
&\ge
\Phi(T_n)-\Phi(T_j)\ge
\frac{1}{2\eta}
\sum_{\ell=j}^{n-1}\|T_{\ell+1}-T_\ell\|^2\ge
\frac{n-j}{2\eta}\|T_n-T_{n-1}\|^2.
\end{aligned}
\]
Summing over $j=1,\ldots,n-1$ and applying \eqref{eq:sum_obj} on the
sum of objective gaps, we obtain
\[
\frac{n(n-1)}{4\eta}\|T_n-T_{n-1}\|^2
\le
\sum_{j=1}^{n-1}\bigl[\Phi^\star-\Phi(T_j)\bigr]
\le
\frac{\|T_0-T^\star\|^2}{2\eta}.
\]
It follows that
\[
\|T_n-T_{n-1}\|
\le
\sqrt{\frac{2}{n(n-1)}}\,\|T_0-T^\star\|
\le
\frac{2\|T_0-T^\star\|}{n},
\qquad n\ge2.
\]
For $n=1$, Fej\'er monotonicity gives
\[
\|T_1-T_0\|
\le
\|T_1-T^\star\|+\|T_0-T^\star\|
\le
2\|T_0-T^\star\|.
\]
Thus, for every $n\ge1$,
\[
\|T_n-T_{n-1}\|
\le
\frac{2\|T_0-T^\star\|}{n}.
\]

Projection optimality in the update defining $T_n$ gives
\[
\langle
T_{n-1}+\eta\nabla\Phi(T_{n-1})-T_n,\,
T'-T_n
\rangle
\le0.
\]
Hence, for every $T'\in\mathcal K$,
\[
\begin{aligned}
\eta\langle\nabla\Phi(T_n),T'-T_n\rangle
&\le
\|
T_n-T_{n-1}
+\eta(\nabla\Phi(T_n)-\nabla\Phi(T_{n-1}))
\|
\|T'-T_n\|\le
2\rho_{\mathcal K}\|T_n-T_{n-1}\|,
\end{aligned}
\]
where the last inequality follows from
\eqref{eq:forward-nonexpansion} and
$\|T'-T_n\|\le2\rho_{\mathcal K}$.
Taking the maximum over $T'\in\mathcal K$ and using the step bound
above yields
\[
\mathcal G(T_n)
\le
\frac{2\rho_{\mathcal K}}{\eta}\|T_n-T_{n-1}\|
\le
\frac{4\rho_{\mathcal K}\|T_0-T^\star\|}{\eta n},
\qquad n\ge1.
\]

Finally, $\Phi(T_n)\to\Phi^\star$ and $\mathcal G(T_n)\to0$ imply
$\bar R(T_n)\to\Phi^\star$.
Since $T_n\to\bar T$ and $J_S$ is continuous,
\[
\bar R(T_n)-J_S(T_n)
\rightarrow
\Phi^\star-J_S(\bar T)
=
\Phi(\bar T)-J_S(\bar T).
\]
Substituting the definitions of the corrected objectives proves
\eqref{eq:limiting-certificate-width}.
\end{proof}

Theorem~\ref{thm:generator-first-order} separates optimization accuracy
from risk certification. Since
$\Phi^\star-\Phi(T_n)\le\mathcal G(T_n)$,
the stopping rule $\mathcal G(T_n)\le\varepsilon_{\rm opt}$ certifies
$\varepsilon_{\rm opt}$-optimality for the corrected problem.
Meanwhile, $J_S(T_n)\le R_S\le\bar R(T_n)$, so the reported interval
provides certified bounds on the original structured worst-case risk.
Any nonzero limiting width equals the residual correction
$\Phi(\bar T)-J_S(\bar T)$.

\subsection{Exact recovery and linear convergence}
\label{sec:exact-generation}

We now characterize when the limiting risk interval closes,
yielding $J_S(\bar T)=R_S=\Phi^\star$.
By Theorem~\ref{thm:generator-first-order}, the limiting width is
$\Phi^\star-J_S(\bar T)$.
In the moment case, the width
vanishes if and only if $\tr[\Omega(B_S-M(\bar T))]=0$, or equivalently
$(B_S-M(\bar T))\Omega=0$: the budget must be saturated on
$\range(\Omega)$, but not necessarily on its orthogonal complement.
In that case $M(\bar T)\Omega=B_S\Omega$, where the right-hand side has
rank $\rank(\Omega)$ because $B_S\succ0$ and the left-hand side has rank
at most $\rank M(\bar T)\le N$. Hence this requires $N\ge\rank(\Omega)$.
For $\Omega\succ0$ in the moment case and $L>0$ in
the Wasserstein case, the interval closes if and only if
\[
M(\bar T)=B_S
\quad\text{(Moment)},
\qquad
\|\bar T\|_{L^2}=\rho
\quad\text{(Wasserstein)},
\]
and the moment condition requires $N\ge k$.
We next give a strict multiplier condition that guarantees exact
recovery and linear convergence.

Let $T^\star$ maximize $\Phi$ over $\mathcal K$, and let $\Lambda^\star$ or
$\lambda^\star$ be an associated restricted-dual multiplier satisfying
\eqref{eq:corrected-matrix-kkt} or \eqref{eq:corrected-scalar-kkt}.
We consider the strict-margin condition
\begin{equation}\label{eq:strict-margin}
\Lambda^\star\succ\Omega
\quad\text{(Moment)},
\qquad
\lambda^\star>\frac L2
\quad\text{(Wasserstein)},
\end{equation}
and measure the margin in the inner product
\eqref{eq:response-inner-product} by
\[
\nu:=
\begin{cases}
\lambda_{\min}\!\left(
B_S^{1/2}(\Lambda^\star-\Omega)B_S^{1/2}
\right), &\text{(Moment)},\\[2mm]
\lambda^\star-\dfrac L2, &\text{(Wasserstein)}.
\end{cases}
\]

\begin{corollary}[Linear convergence under a strict multiplier margin]
\label{cor:strict-generation}
Suppose \eqref{eq:regular} and \eqref{eq:strict-margin} hold.
Then $T^\star$ saturates the corresponding budget and satisfies
$
J_S(T^\star)=R_S=\Phi^\star$.
Moreover, $T^\star$ is the unique maximizer of $\Phi$ over $\mathcal K$, and
$q_{T^\star}$ is the unique optimal conditional distribution for $R_S$. Let $\{T_n\}_{n\ge0}$ be generated by \eqref{eq:projected-generation}
with $T_0\in\mathcal K$ and a constant stepsize $\eta>0$ satisfying
$\eta L_\Phi\le1$. Then, for every $n\ge0$,
\begin{align}
\|T_n-T^\star\|
&\le \vartheta^n\|T_0-T^\star\|,
\label{eq:global-projected-linear}\\
0\le\bar R(T_n)-J_S(T_n)
&\le \kappa\vartheta^n\|T_0-T^\star\|,
\label{eq:global-certificate-rate}
\end{align}
where $\vartheta:=(1+\eta\nu)^{-1}<1$ and
$\kappa:=3L_\Phi \rho_{\mathcal K}+\|\nabla\Phi(T^\star)\|$.
\end{corollary}

\begin{proof}[Proof of Corollary \ref{cor:strict-generation}]
Recall that the norms and inner products are those of
\eqref{eq:response-inner-product}. 
In the moment case, complementary
slackness gives
\[
\tr\!\left[(\Lambda^\star-\Omega)(B_S-M(T^\star))\right]=0.
\]
Since $\Lambda^\star-\Omega\succ0$ and $B_S-M(T^\star)\succeq0$, we obtain
$M(T^\star)=B_S$.
Similarly, in the Wasserstein case,
$(\lambda^\star-L/2)(\rho^2-\|T^\star\|_{L^2}^2)=0$ and $\lambda^\star>L/2$ imply
$\|T^\star\|_{L^2}=\rho$. Thus the correction vanishes, and
\[
J_S(T^\star)\le R_S\le\Phi^\star=\Phi(T^\star)=J_S(T^\star).
\]

The strict margin also gives the following inequality:
\begin{equation}\label{eq:normal-inequality}
\langle\nabla\Phi(T^\star),T^\star-T\rangle
\ge\nu\|T-T^\star\|^2,
\qquad T\in\mathcal K.
\end{equation}
To prove it in the moment case, write $c_i:=\vecop(T_i)$ and
$c_i^\star:=\vecop(T_i^\star)$. 
Stationarity and \eqref{eq:response-gradient}
give
$\vecop((\nabla\Phi(T^\star))_i)=2B_S(\Lambda^\star-\Omega)c_i^\star$.
Using $M(T^\star)=B_S$, we have
\[
\begin{aligned}
\langle\nabla\Phi(T^\star),T^\star-T\rangle
&=2\sum_i p_i c_i^{\star\top}(\Lambda^\star-\Omega)(c_i^\star-c_i)\\
&=\tr\left[(\Lambda^\star-\Omega)M(T-T^\star)\right]+\tr\left[(\Lambda^\star-\Omega)(B_S-M(T))\right]\\
&\ge\nu\tr[B_S^{-1}M(T-T^\star)]
=\nu\|T-T^\star\|^2.
\end{aligned}
\]
The inequality uses $M(T)\preceq B_S$ and
$\Lambda^\star-\Omega\succeq\nu B_S^{-1}$.
In the Wasserstein case, stationarity gives
$\nabla\Phi(T^\star)=2\nu T^\star$. Since $\|T^\star\|=\rho$,
\[
\langle\nabla\Phi(T^\star),T^\star-T\rangle
=\nu\bigl(\|T-T^\star\|^2+\rho^2-\|T\|^2\bigr)
\ge\nu\|T-T^\star\|^2.
\]
This proves \eqref{eq:normal-inequality} for both budgets.

By concavity and \eqref{eq:normal-inequality},
\[
\Phi^\star-\Phi(T)
\ge\langle\nabla\Phi(T^\star),T^\star-T\rangle
\ge\nu\|T-T^\star\|^2,
\qquad T\in\mathcal K.
\]
Hence $T^\star$ is the unique maximizer of $\Phi$.
To prove uniqueness of the worst-case distribution, let $q$ be feasible
with $J_S(q)=R_S$, and let $T_i:=\E_q[C\mid\xi_i]$.
The conditional-mean argument in
Proposition~\ref{prop:exact-shift-recovery} gives $T\in\mathcal K$ and
$J_S(q)\le\Phi(T)$. Therefore
\[
R_S=J_S(q)\le\Phi(T)
\le\Phi^\star-\nu\|T-T^\star\|^2
=R_S-\nu\|T-T^\star\|^2,
\]
so $T=T^\star$. Feasibility of $q$ and budget saturation at $T^\star$ now give
\[
0\le\sum_i p_i\E_q[\|C-T_i^\star\|_F^2\mid\xi_i]
=\E_q\|C\|_F^2-\|T^\star\|_{L^2}^2\le0.
\]
Thus $q(\cdot\mid\xi_i)=\delta_{T_i^\star}$ for every $i$, proving
$q=q_{T^\star}$.

For the convergence rate, projection optimality gives
\[
\langle T_n+\eta\nabla\Phi(T_n)-T_{n+1},\,
T^\star-T_{n+1}\rangle\le0.
\]
Combining this with \eqref{eq:normal-inequality} at $T=T_{n+1}$ yields
\[
\begin{aligned}
(1+\eta\nu)\|T_{n+1}-T^\star\|^2
&\le\langle
(T_n+\eta\nabla\Phi(T_n))-(T^\star+\eta\nabla\Phi(T^\star)),\,
T_{n+1}-T^\star\rangle\\
&\le\|T_n-T^\star\|\,\|T_{n+1}-T^\star\|,
\end{aligned}
\]
where the last inequality follows from Cauchy--Schwarz and
\eqref{eq:forward-nonexpansion}. Hence
$\|T_{n+1}-T^\star\|\le\vartheta\|T_n-T^\star\|$, which proves
\eqref{eq:global-projected-linear}.

Finally, budget saturation at $T^\star$ gives, for every $T\in\mathcal K$,
\[
0\le\Phi(T)-J_S(T)
=\sum_ip_i\bigl[c_i^{\star\top}\Omega c_i^\star-c_i^\top\Omega c_i\bigr]
=\sum_ip_i(c_i^\star-c_i)^\top\Omega(c_i^\star+c_i)
\]
in the moment case, and the same identity with $\Omega=\frac L2I_k$ in the
Wasserstein case. By Cauchy--Schwarz and
$2\sum_ip_i\vecop(T_i')^\top\Omega\vecop(T_i')\le L_\Phi\|T'\|^2$, which
holds for every $T'$ by the definition of $L_\Phi$,
\[
0\le\Phi(T)-J_S(T)
\le\frac{L_\Phi}{2}\|T-T^\star\|\,\|T+T^\star\|
\le L_\Phi \rho_{\mathcal K}\|T-T^\star\|.
\]
Also, first-order optimality at $T^\star$ gives, for $T'\in\mathcal K$,
\[
\begin{aligned}
\langle\nabla\Phi(T),T'-T\rangle
&\le\langle\nabla\Phi(T)-\nabla\Phi(T^\star),T'-T\rangle
+\langle\nabla\Phi(T^\star),T^\star-T\rangle\\
&\le\bigl(2L_\Phi \rho_{\mathcal K}+\|\nabla\Phi(T^\star)\|\bigr)
\|T-T^\star\|.
\end{aligned}
\]
Taking the maximum over $T'\in\mathcal K$ and adding the correction
bound, we obtain
\[
\bar R(T)-J_S(T)
=\Phi(T)-J_S(T)+\mathcal G(T)
\le \kappa\|T-T^\star\|.
\]
Applying \eqref{eq:global-projected-linear} proves
\eqref{eq:global-certificate-rate}.
\end{proof}

The linear rate follows from \eqref{eq:normal-inequality}, which
is implied by the strict multiplier margin and budget saturation;
strong concavity of $\Phi$ is not required.
Moreover,
\[
 0\le\mathcal G(T_n)
 \le\bar R(T_n)-J_S(T_n)
 \le\kappa\vartheta^n\|T_0-T^\star\|.
\]
For a fixed problem and stepsize,
Algorithm~\ref{alg:regular-generation} therefore reaches
$\mathcal G(T_n)\le\varepsilon_{\rm opt}$ in
$\mathcal O(\log(1/\varepsilon_{\rm opt}))$ iterations as
$\varepsilon_{\rm opt}\downarrow0$.
For the full support $S=I_d$, the Wasserstein strict-margin condition is
the calibration condition \eqref{eq:calibrated-response}: stationarity in
\eqref{eq:corrected-scalar-kkt} with $\lambda^\star>L/2$ and budget
saturation give $Z_i^*=T_i^\star$. Hence
Algorithm~\ref{alg:regular-generation} computes the
calibrated responses of Section~\ref{sec:w-structure} at a linear rate.

\begin{remark}[Verifying strict multiplier margin]
\label{rem:margin-check}
For the Wasserstein budget, a sufficient condition can be checked before
optimization as in Lemma~\ref{lem:calibration}. Let $g(0):=(g_i(0))_{i=1}^N$. If $L>0$ and
\[
0<\rho<\frac{\|g(0)\|_{L^2}}{2L},
\]
then the strict margin holds with
$\nu\ge\|g(0)\|_{L^2}/(2\rho)-L$. Indeed, since $\nabla\Phi_{\rm W}$ is
$2L$-Lipschitz and $\nabla\Phi_{\rm W}(0)=g(0)$,
\[
\|\nabla\Phi_{\rm W}(T)\|_{L^2}
\ge
\|g(0)\|_{L^2}-2L\|T\|_{L^2}
\ge
\|g(0)\|_{L^2}-2L\rho>0
\quad\text{for every}\ T\in\mathcal K_{\rm W}.
\]
Thus no maximizer can lie in the interior of $\mathcal K_{\rm W}$.
At a maximizer $T^\star$, we therefore have $\|T^\star\|_{L^2}=\rho$.
The KKT condition
$\nabla\Phi_{\rm W}(T^\star)=2(\lambda^\star-L/2)T^\star$, with
$\lambda^\star\ge L/2$, gives
\[
\nu=\lambda^\star-\frac L2
=\frac{\|\nabla\Phi_{\rm W}(T^\star)\|_{L^2}}{2\rho}
\ge
\frac{\|g(0)\|_{L^2}-2L\rho}{2\rho}>0.
\]

For either budget, the associated restricted-dual multiplier can also
be recovered at a budget-saturating maximizer $T^\star$ of $\Phi$:
\[
\begin{aligned}
\Lambda^\star
&=
\frac12\left(
\sum_i p_i\vecop(g_i(T_i^\star))\vecop(T_i^\star)^\top
\right)B_S^{-1},
&&\text{(Moment)},\\[1mm]
\lambda^\star
&=
\frac{1}{2\rho^2}
\sum_i p_i\langle g_i(T_i^\star),T_i^\star\rangle_F,
&&\text{(Wasserstein)}.
\end{aligned}
\]
These identities follow from the stationarity equations and
$M(T^\star)=B_S$ or $\|T^\star\|_{L^2}^2=\rho^2$, respectively.
The recovered multiplier can then be tested against
\eqref{eq:strict-margin}. At a finite iterate, these formulas give
estimates of the multiplier.
\end{remark}

\section{Numerical Experiments}\label{sec:experiments}

In this section, we evaluate the proposed support constructions, excluded-risk bounds, and worst-case generation algorithm on a spatial precipitation forecasting task
using precipitation maps from the ERA5 reanalysis dataset
\citep{store2024era5}.

\paragraph{Data and model.}
The data consist of $8{,}760$ hourly precipitation maps. Each map is an $80\times80$ grid, so
$m_1=m_2=80$ and $d=6{,}400$. We partition the grid into 16 contiguous
$20\times20$ blocks and predict their average precipitation six hours
ahead. 
Let $U_0,V_0\in\R^{80\times4}$ have column-orthonormal
block-averaging columns, with each column equal to $1/\sqrt{20}$
on one of four consecutive groups of 20 rows or columns and zero
elsewhere. Consequently, each entry of
$U_0^\top(X-\mu)V_0$ equals $1/20$ times the centered sum over the
corresponding $20\times20$ block, where $\mu\in\R^{80\times80}$ is the cellwise mean map over the fitting period.
For each of the 16 regions, we use its current average precipitation as a baseline and fit a ridge-regression correction using the current averages of all 16 regions as predictors. We use the first half of the chronological record for model fitting and its last quarter to select the ridge penalty. At both the validation and evaluation boundaries, we remove training pairs whose six-hour-ahead targets reach or cross the next period. After selecting the ridge penalty, we fit the final predictor on all 4,371 fitting observations. This yields an affine predictor with parameters $(A,b)\in\R^{16\times16}\times\R^{16}$. The second half of the record is reserved for evaluation. Sampling it every six hours gives $N=730$ observations $\xi_i=(X_i,y_i)$ with equal weights $p_i=1/N$.

For a perturbation $Z\in\R^{80\times80}$ of $X_i$, define the forecast
residual
\[
 e_i(Z):=A\vecop(U_0^\top(X_i+Z-\mu)V_0)+b-y_i,
\]
and we use the quadratic loss
\begin{equation}\label{eq:exp-quadratic-loss}
 f_i(Z)=\frac12\|e_i(Z)\|^2.
\end{equation}
It has the feature structure \eqref{eq:ridge} with
$(r_0,s_0)=(4,4)$, $S_0=V_0\otimes U_0\in\R^{6400\times16}$, and
$
 h_i(C_0)=\frac12\|A\vecop(C_0)+e_i(0)\|^2$.
The common feature Hessian $A^\top A$ is positive definite for the fitted
predictor. Thus, \eqref{eq:ridge} and \eqref{eq:regular} hold with
$L=\lambda_{\max}(A^\top A)$ and
$\alpha=\lambda_{\min}(A^\top A)>0$.

\paragraph{Budgets and risk measures.}

For the moment budget, we fix the spatial geometry $B_{c}\otimes B_{r}$ with
\[
 (B_{c})_{jk}=(B_{r})_{jk}
 =\exp\!\left[-\frac{(j-k)^2}{2\cdot10^2}\right]
 +0.05\,\delta_{jk},
\]
which allows stronger correlation between perturbations at nearby grid
locations. We then calibrate the overall moment scale
$\sigma_{\rm cal}>0$ and the Wasserstein radius $\rho_{\rm cal}>0$ on
the fitting data so that, under each budget, the unrestricted worst-case
quadratic risk equals twice the nominal risk. We set
$B_{\rm cal}:=\sigma_{\rm cal}^2(B_{c}\otimes B_{r})$ and use
$B_{\rm cal}$ and $\rho_{\rm cal}$ as reference budgets throughout
the experiments.

For each loss and budget, write $R_0:=\sum_i p_i f_i(0)$ for the nominal
risk and $R$ and $R_S$ for the unrestricted and support-restricted
worst-case risks. We report both the excluded risk $R-R_S$ and the
excluded risk fraction $(R-R_S)/(R-R_0)$. The latter measures the share
of the worst-case risk increase lost through the support restriction
and equals zero for a lossless support.

\begin{figure}[tbp]
\centering
\includegraphics[width=0.8\textwidth]{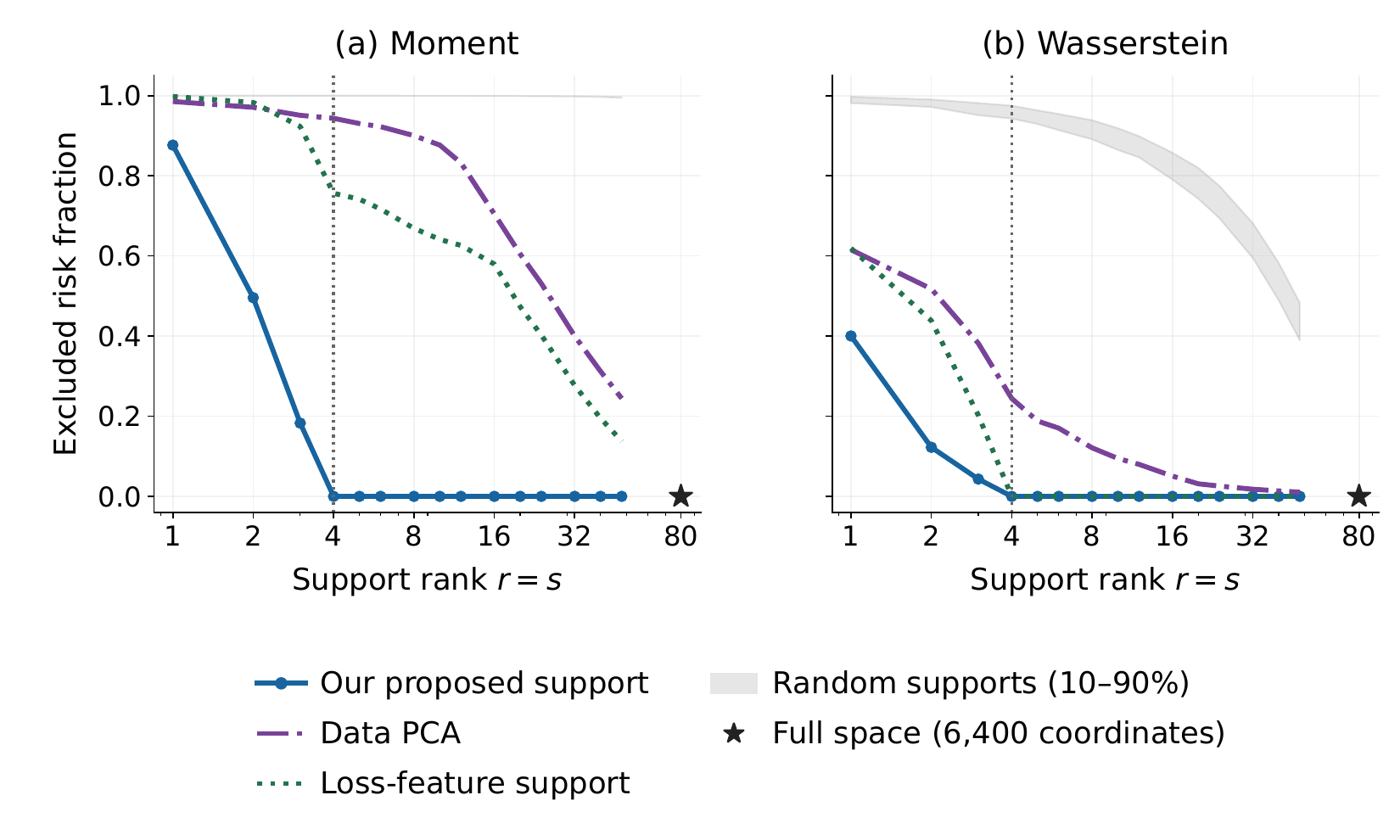}
\caption{Excluded risk fraction $(R-R_S)/(R-R_0)$ versus the common row and column rank $r=s$ under the moment and Wasserstein budgets. The dotted line marks the loss-feature rank $r=s=4$, and the star denotes the unrestricted $80\times80$ perturbation space. The shaded band shows the 10th--90th percentiles over 200 random column-orthonormal supports at each rank.}
\label{fig:exp-support}
\end{figure}

\paragraph{Lossless supports.}
Figure~\ref{fig:exp-support} compares the proposed supports from
Corollaries~\ref{cor:least} and~\ref{cor:w-least} with three baselines over
equal row and column ranks $r=s$. Data PCA uses the leading row and column
subspaces of the centered nominal maps. The loss-feature baseline uses $U_0,V_0$ at $r=s=4$ and is extended
by data-PCA directions at larger ranks. For the loss-feature baseline at ranks below four, we report the median over 20 random pairs of subspaces within the feature spaces. For the random-support baseline, we sample 200 independent pairs of column-orthonormal factors at each rank and show the 10th--90th percentiles.

The proposed supports preserve the unrestricted risk to numerical precision
at $r=s=4$ under both budgets, reducing the number of perturbation
coordinates per observation from $6{,}400$ to 16. At the same rank, data PCA
excludes $94.3\%$ of the worst-case risk increase under the moment budget
and $24.4\%$ under the Wasserstein budget, showing that directions of
dominant nominal variation need not preserve worst-case risk. The random-support band shows that generic low-rank orientations fail to preserve a substantial fraction of the worst-case risk. The loss-feature baseline further highlights the role
of budget geometry. It is lossless under the Wasserstein budget because any
perturbation can be projected onto the loss feature space without changing
the loss or increasing its Frobenius norm, but it excludes $75.6\%$ of the
risk increase under the moment budget. In the latter case, losslessness
depends on the budget-transformed directions $\range(BS_0)$ rather than
$\range(S_0)$ alone, as characterized in Theorem~\ref{thm:ridge} and
Remark~\ref{rem:loss-features}.

\begin{figure}[tbp]
\centering
\begin{subfigure}{0.24\textwidth}
\parbox[t][1.35\linewidth][t]{\linewidth}{\vspace{0pt}%
\includegraphics[width=\linewidth]{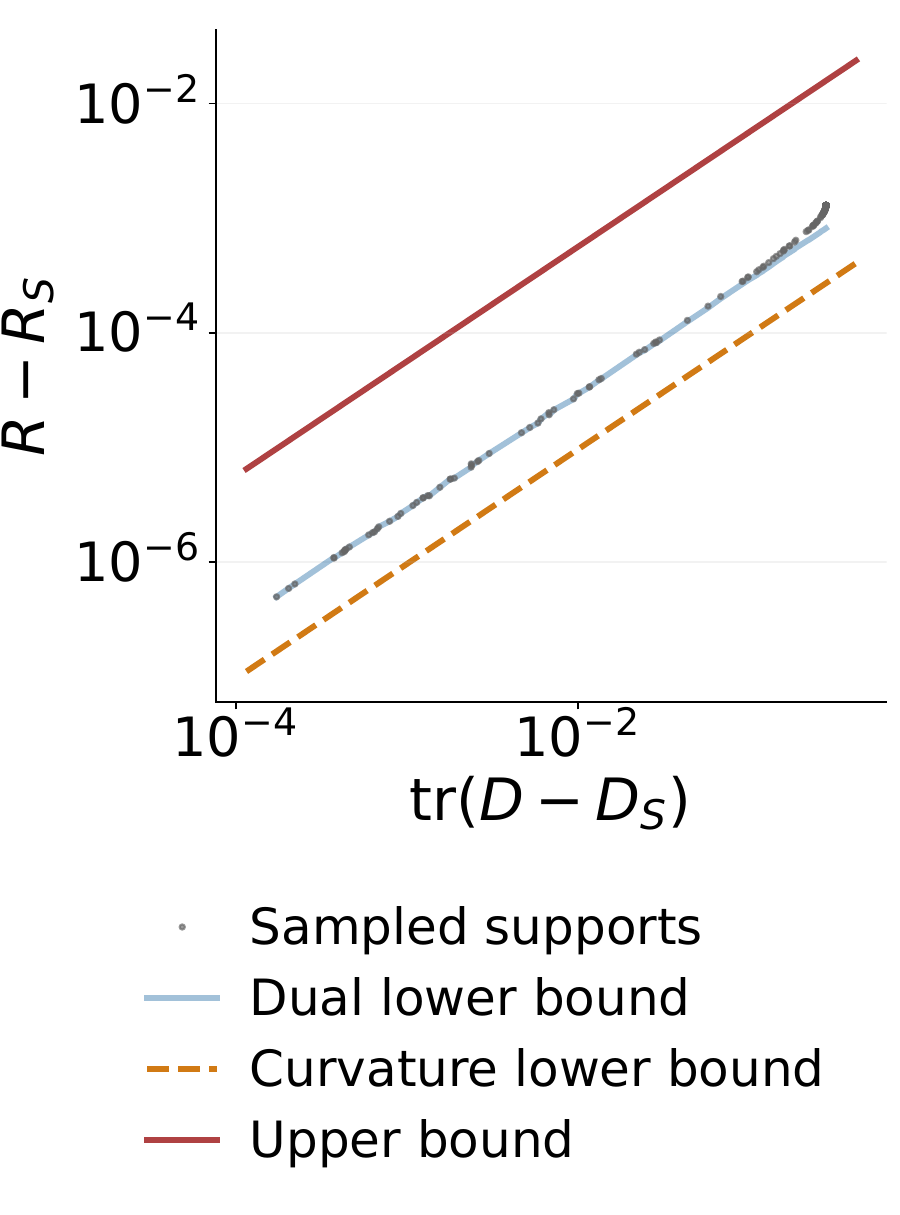}}
\par\smallskip{\small (a)}
\end{subfigure}\hfill%
\begin{subfigure}{0.24\textwidth}
\parbox[t][1.35\linewidth][t]{\linewidth}{\vspace{0pt}%
\includegraphics[width=\linewidth]{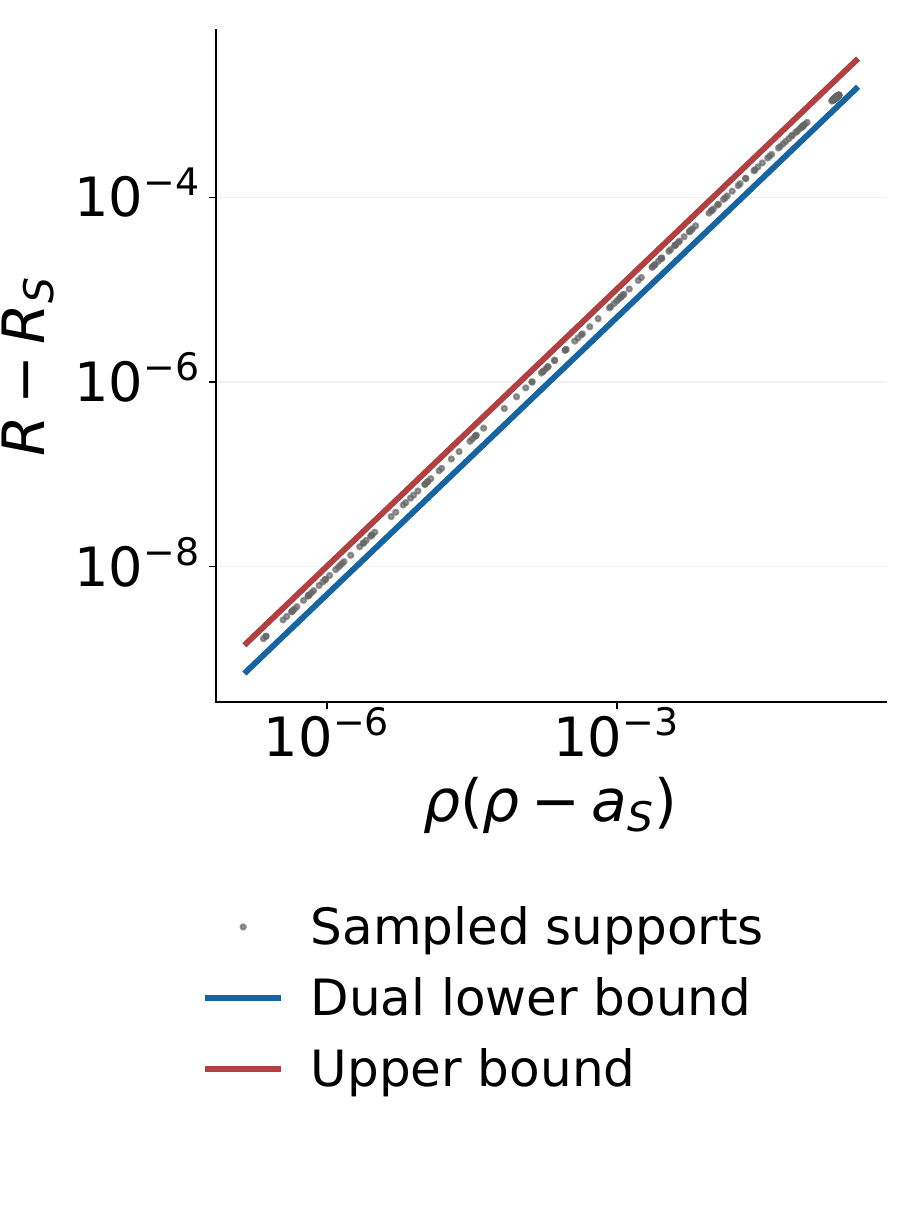}}
\par\smallskip{\small (b)}
\end{subfigure}\hfill%
\begin{subfigure}{0.24\textwidth}
\parbox[t][1.35\linewidth][t]{\linewidth}{\vspace{0pt}%
\includegraphics[width=\linewidth]{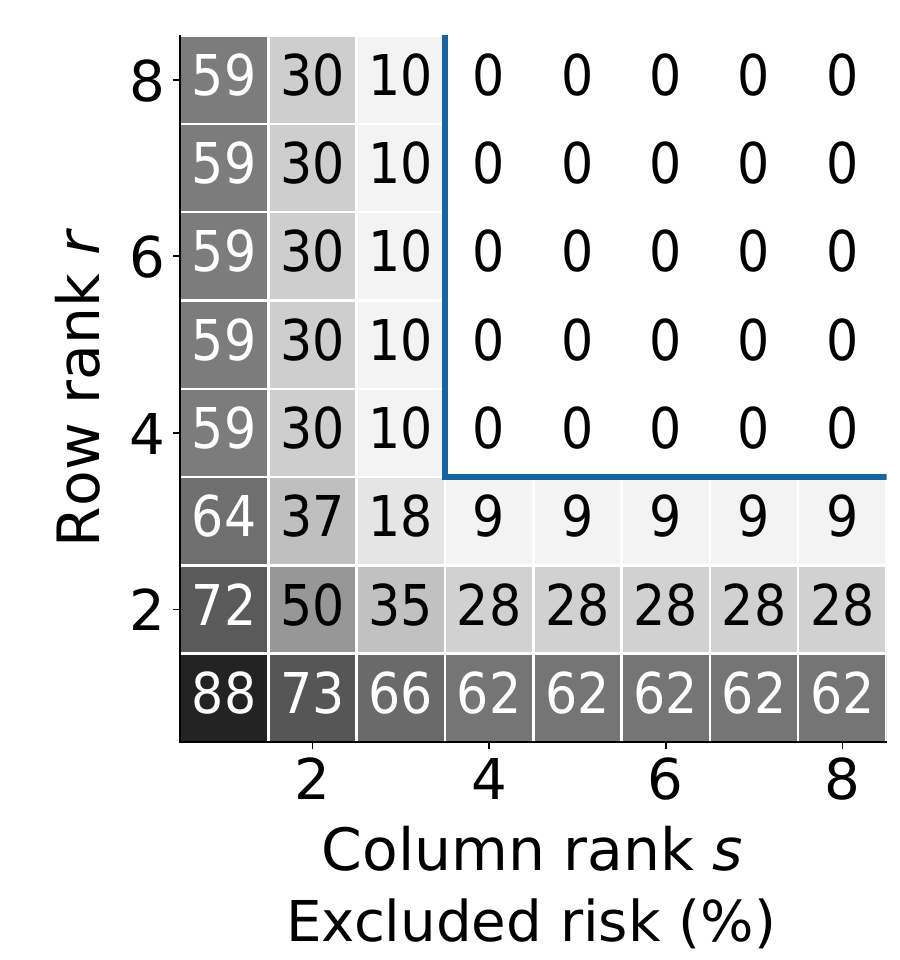}}
\par\smallskip{\small (c)}
\end{subfigure}\hfill%
\begin{subfigure}{0.24\textwidth}
\parbox[t][1.35\linewidth][t]{\linewidth}{\vspace{0pt}%
\includegraphics[width=\linewidth]{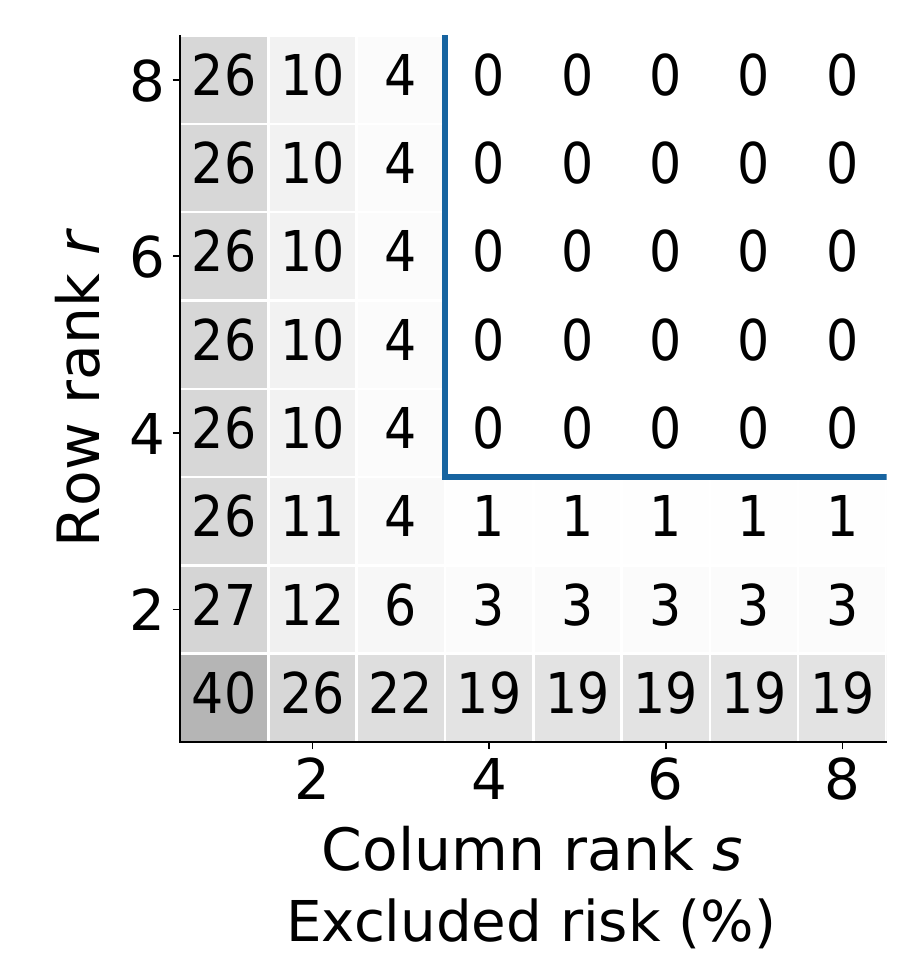}}
\par\smallskip{\small (d)}
\end{subfigure}
\caption{(a,b) Excluded risk $R-R_S$ for 300 sampled supports under each budget, plotted against $\tr(D-D_S)$ for the moment budget and $\rho(\rho-a_S)$ for the Wasserstein budget, together with the bounds of Theorems~\ref{thm:ridge} and~\ref{prop:energy}. In panel~(a), gray dots show the excluded risk, the light-blue line joins the support-specific dual lower bounds $\tr[\Lambda_0(D-D_S)]$, and the dashed line shows the curvature bound $(\alpha/2)\tr(D-D_S)$. (c,d) Excluded risk percentage $100(R-R_S)/(R-R_0)$ over rank pairs $(r,s)$ for the proposed supports, rounded to integers. The lines mark the predicted minimal lossless ranks $r=4$ and $s=4$.}
\label{fig:exp-bounds}
\end{figure}

\paragraph{Excluded risk bounds and minimal lossless dimensions.}
Figure~\ref{fig:exp-bounds}(a,b) tests the quantitative bounds on the
excluded risk in Theorems~\ref{thm:ridge} and~\ref{prop:energy} using
300 supports of varying ranks and orientations under each budget. The
horizontal coordinates are $\tr(D-D_S)$ for the moment budget and
$\rho(\rho-a_S)$ for the Wasserstein budget. All computed gaps lie between
their corresponding upper and lower bounds. Under the moment budget, the
dual lower bound $\tr[\Lambda_0(D-D_S)]$ closely tracks the excluded risk
and is substantially tighter than the curvature bound
$(\alpha/2)\tr(D-D_S)$. The multiplier $\Lambda_0$ is computed once at the
full feature budget $D$ and then used for every support. Unlike the
curvature bound, which prices all lost directions uniformly by
$\alpha/2$, the dual bound weights them according to the full-problem
multiplier $\Lambda_0$. Panels~(c,d) scan
$(r,s)\in\{1,\ldots,8\}^2$ and recover the predicted minimal lossless
dimensions $r_*=s_*=4$, corresponding to the region $r\ge4$ and $s\ge4$.

\paragraph{Computational cost.}
We next examine the computational benefit of the low-rank supports.
Table~\ref{tab:exp-cost} compares Algorithm~\ref{alg:regular-generation}
on the proposed supports with $r=s\in\{1,4,8,16\}$ and on the unrestricted
$80\times80$ space. Let
\[
 A_S:=AS_0^\top S,
 \qquad
 H_S:=A_S^\top A_S
 =S^\top S_0A^\top AS_0^\top S .
\]
For the quadratic loss, $H_S$ is the constant Hessian of the restricted
loss, and we use the exact moment curvature correction
$\Omega=H_S/2$. Each run starts from $T_0=0$ and stops when
$
 (\bar R(T_n)-J_S(T_n))/{(R_S-R_0)}\le10^{-6}$.
All computations
are performed in the $k=rs$ dimensional core using the corresponding
feature Hessian and induced budget.

At the minimal lossless rank $r=s=4$, the algorithm reaches the
$10^{-6}$ tolerance in 6 Wasserstein iterations and 103 moment iterations,
compared with 6 and 111 iterations in the full space. The resulting
speedup therefore comes primarily from the lower per-iteration cost.
Total solver time decreases from $0.669$ s to $0.956$ ms under the
Wasserstein budget and from $67.0$ s to $0.129$ s under the moment
budget, corresponding to speedups of approximately $700$ and $520$, respectively, while preserving the unrestricted risk to numerical precision. The rank-one results illustrate the further computational
savings available when some loss of worst-case risk is allowed, as
quantified in Figure~\ref{fig:exp-support}.

\begin{table}[htbp]
\centering
\caption{Computational cost for the quadratic loss at relative certificate
tolerance $10^{-6}$. Times are medians over three runs on an Intel Xeon
Gold 6226 processor with four threads and $N=730$. Total time includes all
algorithm iterations and excludes support construction, data loading, and
setup. Here $K$ denotes the iteration count, and time per iteration is the
total time divided by $K$. The bold row corresponds to the minimal lossless
support.}
\label{tab:exp-cost}
\small
\begin{tabular}{@{}lccccccc@{}}
\toprule
 & & \multicolumn{3}{c}{Wasserstein}
 & \multicolumn{3}{c}{Moment}\\
\cmidrule(lr){3-5}\cmidrule(lr){6-8}
$r=s$ & $k=rs$
& per iter. (ms) & $K$ & total (s)
& per iter. (ms) & $K$ & total (s)\\
\midrule
$1$ & $1$
& $0.134$ & $2$ & $0.000267$
& $0.406$ & $2$ & $0.000811$\\
$\mathbf{4}$ & $\mathbf{16}$
& $\mathbf{0.159}$ & $\mathbf{6}$ & $\mathbf{0.000956}$
& $\mathbf{1.25}$ & $\mathbf{103}$ & $\mathbf{0.129}$\\
$8$ & $64$
& $0.488$ & $6$ & $0.00293$
& $6.37$ & $103$ & $0.656$\\
$16$ & $256$
& $2.05$ & $6$ & $0.0123$
& $34.7$ & $103$ & $3.58$\\
\midrule
Full ($80$) & $6{,}400$
& $112$ & $6$ & $0.669$
& $604$ & $111$ & $67.0$\\
\bottomrule
\end{tabular}
\end{table}

\paragraph{Worst-case generation.}
We next evaluate Algorithm~\ref{alg:regular-generation} on a data-driven
support with $r=s=4$ learned from the fitting maps by
CoreFlow~\citep{wu2026coreflow}, providing a fixed support chosen
independently of the loss and uncertainty budget. We vary the budget scale by setting
\[
 \rho=\gamma\rho_{\rm cal},
 \qquad
 B=\gamma^2B_{\rm cal},
\]
and use the same calibrated reference budgets for both the quadratic and
Cauchy losses without recalibration. All runs start from $T_0=0$.

\begin{figure}[h]
\centering
\begin{subfigure}{0.24\textwidth}
\includegraphics[width=\linewidth]{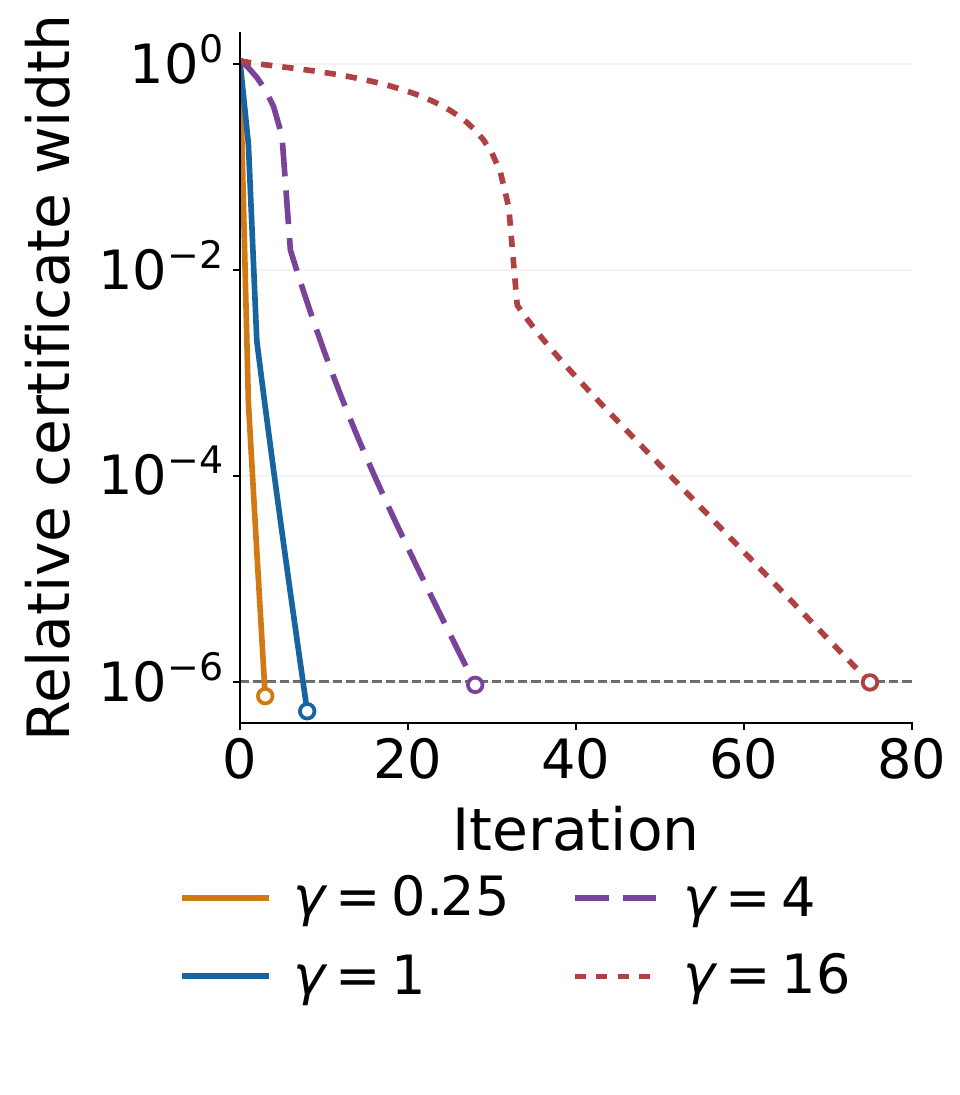}
\par\smallskip{\small (a)}
\end{subfigure}\hfill%
\begin{subfigure}{0.24\textwidth}
\includegraphics[width=\linewidth]{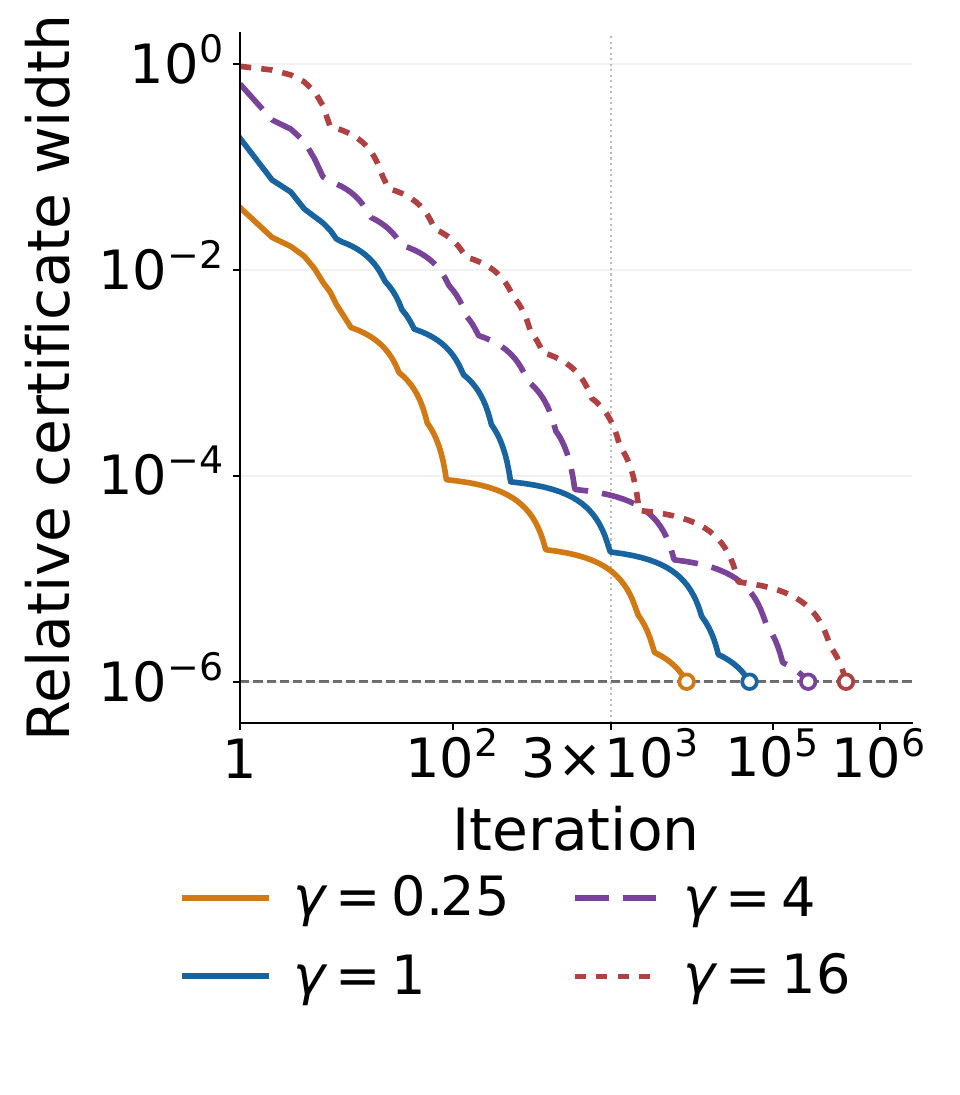}
\par\smallskip{\small (b)}
\end{subfigure}\hfill%
\begin{subfigure}{0.24\textwidth}
\includegraphics[width=\linewidth]{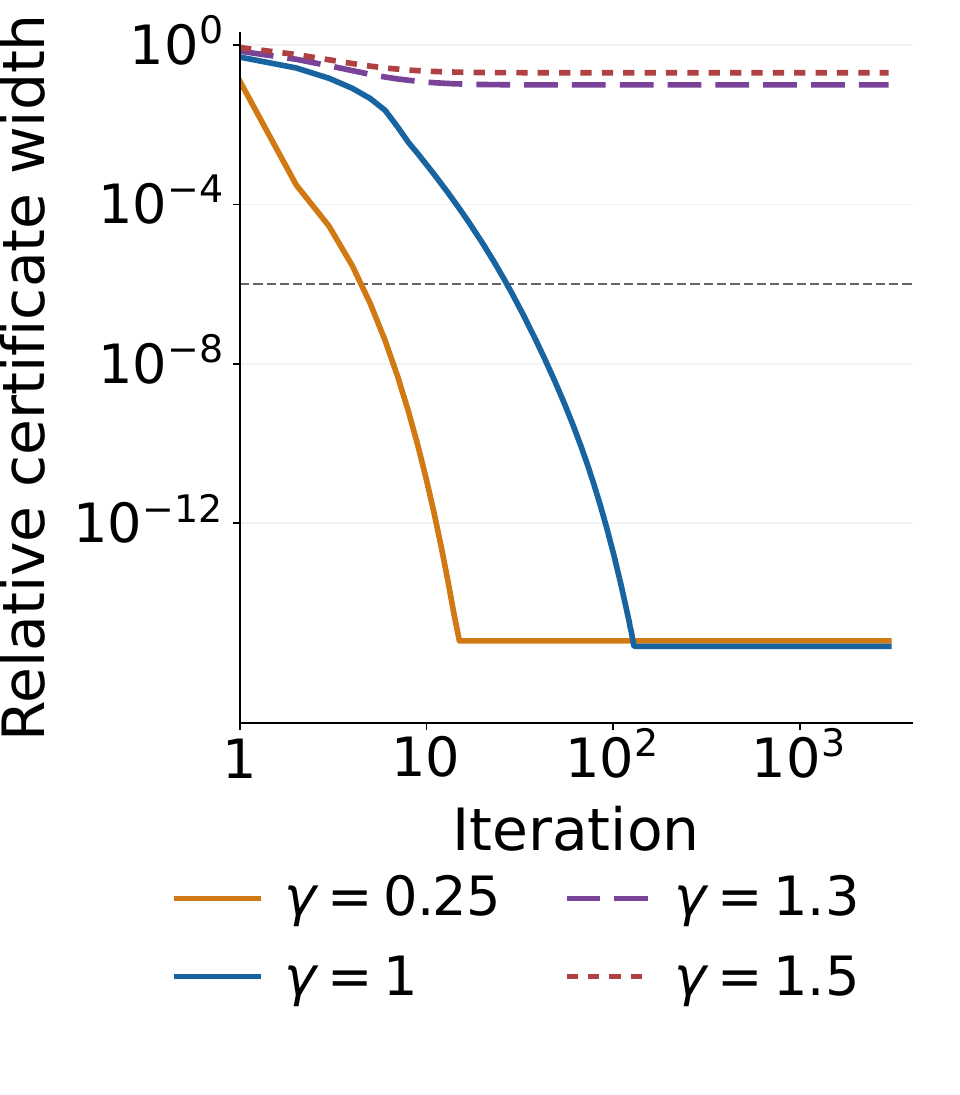}
\par\smallskip{\small (c)}
\end{subfigure}\hfill%
\begin{subfigure}{0.24\textwidth}
\includegraphics[width=\linewidth]{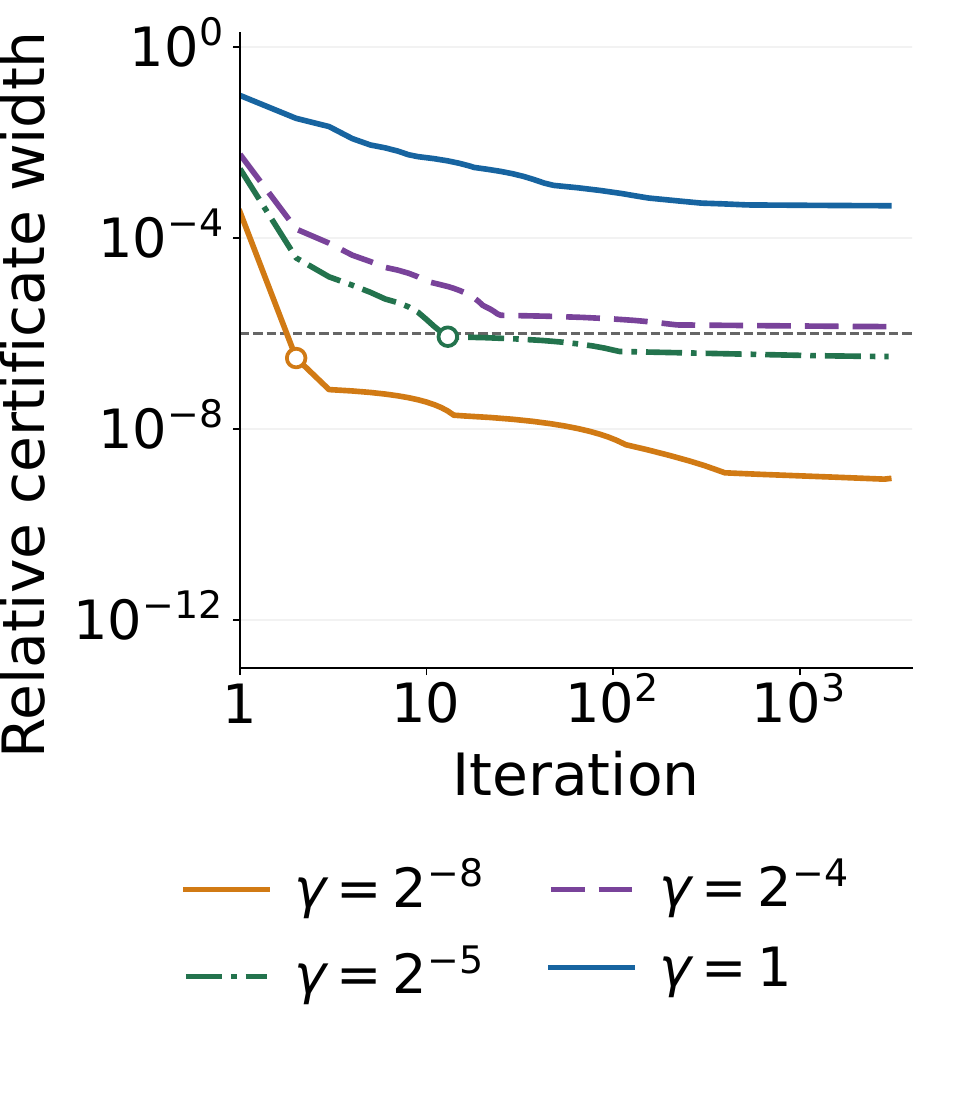}
\par\smallskip{\small (d)}
\end{subfigure}
\caption{Certificate trajectories on the CoreFlow support with $r=s=4$.
Panels~(a,b) use the quadratic loss under the Wasserstein and moment
budgets, respectively, while panels~(c,d) use the Cauchy loss under the
Wasserstein and moment budgets. The plotted widths are
$(\bar R(T_n)-J_S(T_n))/(R_S-R_0)$ in panels~(a,b) and
$(\bar R(T_n)-J_S(T_n))/J_S(T_{3000})$ in panels~(c,d). The dashed lines mark the $10^{-6}$ relative tolerance; circles in
panels~(a,b,d) mark its first crossing. In panels~(c,d), the horizontal
coordinate is $n+1$, with $T_0$ plotted at 1.
For logarithmic display, absolute certificate widths below $10^{-18}$
are plotted at $10^{-18}$ before normalization.}
\label{fig:exp-certificates}
\end{figure}

\emph{Quadratic loss.}
Figure~\ref{fig:exp-certificates}(a,b) considers
\eqref{eq:exp-quadratic-loss} at
$\gamma\in\{0.25,1,4,16\}$. Under the moment budget, we use the exact
curvature correction $\Omega=H_S/2$. The constant stepsizes are
$\eta=(2L_S)^{-1}$ for Wasserstein and
$\eta=[2L_S\lambda_{\max}(B_S)]^{-1}$ for the moment budget, where $L_S=\lambda_{\max}(H_S)$. Since
$0\preceq H_S\preceq L_SI_k$, these choices satisfy the stepsize conditions
of Theorem~\ref{thm:generator-first-order}. The moment correction removes
the quadratic term, leaving the affine objective
\[
 \Phi_{\rm M}(T)
 =R_0+\frac12\tr(H_SB_S)
 +\sum_i p_i\bigl(A_S^\top e_i(0)\bigr)^\top\vecop(T_i).
\]
All 8 configurations reach the $10^{-6}$ relative certificate
tolerance. The 4 Wasserstein trajectories cross the threshold within
75 iterations. The moment trajectories converge more slowly on the learned support, with
the $\gamma=1$ run reaching the threshold at iteration $59{,}950$. The
linear convergence guarantee of Corollary~\ref{cor:strict-generation}
holds at every scale, with a contraction factor that weakens as the budget
scale increases, in line with the longer trajectories at larger $\gamma$.

\emph{Cauchy loss.}
Figure~\ref{fig:exp-certificates}(c,d) keeps the predictor and learned
support fixed and replaces the quadratic loss by
\[
 f_i(Z)=\ell_c(e_i(Z)),
 \qquad
 \ell_c(e)=\frac{c^2}{2}\sum_{j=1}^{16}
 \log\!\left(1+\frac{e_j^2}{c^2}\right),
\]
which downweights large regional residuals. We set $c=2.385\hat\sigma$, where
$\hat\sigma=1.4826\,\mathrm{MAD}$ is the Gaussian-consistent MAD scale
estimate of the fitting residuals. The factor $2.385$ is the standard
Cauchy tuning constant corresponding to approximately $95\%$ efficiency
under Gaussian errors. For
$
 \varphi_c(t)=\frac{c^2}{2}\log\!\left(1+\frac{t^2}{c^2}\right)$,
we have
$
 \varphi_c''(t)
 =\frac{1-t^2/c^2}{(1+t^2/c^2)^2}
 \in[-1/8,1]$.
Consequently, $-\frac18H_S\preceq\nabla^2 f_i^S(c)\preceq H_S$ for every $i$ and $c$, so
$L_S=\lambda_{\max}(H_S)$ is a valid curvature bound. We use this bound
in the Wasserstein correction and stepsize, and take
$\Omega=(L_S/2)I_k$ under the moment budget. Each run is carried out for
$3{,}000$ iterations.

Under the Wasserstein budget, the certificate closes to machine precision
at $\gamma\in\{0.25,1\}$ and stalls at a positive width at
$\gamma\in\{1.3,1.5\}$. Under the moment budget, the two smallest scales
reach the $10^{-6}$ tolerance within a few iterations, while the widths at
$\gamma=2^{-4}$ and $\gamma=1$ level off above it. The decomposition
$\bar R(T_n)-J_S(T_n)=\mathcal G(T_n)+[\Phi(T_n)-J_S(T_n)]$ identifies the
source of these residual widths: for the moment run at $\gamma=1$,
$\mathcal G(T_{3000})/J_S(T_{3000})=3.63\times10^{-6}$, while the
correction $\tr[\Omega(B_S-M(T_{3000}))]$ accounts for $99.2\%$ of the
final width.

For the Cauchy loss, logarithmic growth is dominated by the quadratic
penalty in the isotropic corrections used here, so the corrected objective
has a finite unconstrained maximizer, which does not depend on the budget
scale. Under the Wasserstein budget, the certificate closes while this maximizer
lies outside the ball and has a positive limiting width once it lies
inside; under the moment budget, the width is positive as soon as the
budget is slack in some direction. In both cases the iterates still
converge by Theorem~\ref{thm:generator-first-order} and the optimization
gap vanishes, so the remaining correction is exactly the limiting
certificate width in \eqref{eq:limiting-certificate-width}. 
As the budget approaches this point, the margin $\nu$ shrinks to zero and
the linear rate of Corollary~\ref{cor:strict-generation} slows, as seen
for $\gamma=1$ in Figure~\ref{fig:exp-certificates}(c). Beyond it, only
the sublinear guarantee of Theorem~\ref{thm:generator-first-order}
applies. Thus, increasing the budget can change a closed certificate into
one with a positive limiting width, even though the corrected problem is
solved to optimality.

\section{Closing Remarks}\label{sec:conclusion}

We studied the risk cost of restricting matrix perturbations to shared
row and column spaces under moment and Wasserstein budgets. For a fixed
decision, the analysis bounds the excluded worst-case risk of any
prescribed support, characterizes the lossless supports, and identifies
their minimal dimensions. The lossless directions are determined jointly
by the loss and the budget instead of the dominant variation of the
data, which separates them from principal components. For a fixed
support, the curvature-corrected formulation yields feasible
projected gradient iterates with computable risk certificates. The
certificate closes with linear convergence under a strict multiplier
margin and otherwise converges to the price of the unused budget.
Together, these results separate the risk excluded by the support from
the error in generating a worst-case distribution within that support.

\section*{Acknowledgment}

This work is partially supported by NSF CMMI-2112533 and the Coca-Cola Foundation.

\bibliographystyle{plainnat}
\bibliography{reference}

@article{cook2007dimension,
  author        = {Cook, R. Dennis},
  title         = {{Fisher} Lecture: Dimension Reduction in Regression},
  journal       = {Statistical Science},
  volume        = {22},
  number        = {1},
  pages         = {1--26},
  year          = {2007},
}

@article{li2010dimension,
  author        = {Li, Bing and Kim, Min Kyung and Altman, Naomi},
  title         = {On Dimension Folding of Matrix- or Array-Valued Statistical Objects},
  journal       = {The Annals of Statistics},
  volume        = {38},
  number        = {2},
  pages         = {1094--1121},
  year          = {2010},
}

@misc{store2024era5,
  author        = {Hersbach, H. and Bell, B. and Berrisford, P. and Biavati, G. and Hor{\'a}nyi, A. and Mu{\~n}oz Sabater, J. and Nicolas, J. and Peubey, C. and Radu, R. and Rozum, I. and Schepers, D. and Simmons, A. and Soci, C. and Dee, D. and Th{\'e}paut, J.-N.},
  title         = {{ERA5} hourly data on single levels from 1940 to present},
  year          = {2018},
  howpublished  = {Copernicus Climate Change Service (C3S) Climate Data Store (CDS)},
}

@article{bertsimasSim2004,
  author        = {Dimitris Bertsimas and Melvyn Sim},
  title         = {The Price of Robustness},
  journal       = {Operations Research},
  volume        = {52},
  number        = {1},
  pages         = {35--53},
  year          = {2004},
}

@article{blanchetMurthy2019,
  author        = {Jose Blanchet and Karthyek Murthy},
  title         = {Quantifying distributional model risk via optimal transport},
  journal       = {Mathematics of Operations Research},
  volume        = {44},
  number        = {2},
  pages         = {565--600},
  year          = {2019},
}

@misc{cheng2025,
  author        = {Xiuyuan Cheng and Yao Xie and Linglingzhi Zhu and Yunqin Zhu},
  title         = {Worst-case generation via minimax optimization in {Wasserstein} space},
  year          = {2025},
  howpublished  = {{arXiv} preprint arXiv:2512.08176},
  eprint        = {2512.08176},
  archivePrefix = {arXiv},
}

@article{delage2010,
  author        = {Erick Delage and Yinyu Ye},
  title         = {Distributionally robust optimization under moment uncertainty with application to data-driven problems},
  journal       = {Operations Research},
  volume        = {58},
  number        = {3},
  pages         = {595--612},
  year          = {2010},
}

@article{gaoKleywegt2023,
  author        = {Rui Gao and Anton~J. Kleywegt},
  title         = {Distributionally robust stochastic optimization with {Wasserstein} distance},
  journal       = {Mathematics of Operations Research},
  volume        = {48},
  number        = {2},
  pages         = {603--655},
  year          = {2023},
}

@article{goldfarbIyengar2003,
  author        = {Donald Goldfarb and Garud Iyengar},
  title         = {Robust Portfolio Selection Problems},
  journal       = {Mathematics of Operations Research},
  volume        = {28},
  number        = {1},
  pages         = {1--38},
  year          = {2003},
}

@article{jiang2026optimized,
  author        = {Jiang, Shiyi and Cheng, Jianqiang and Pan, Kai and Shen, Zuo-Jun Max},
  title         = {Optimized dimensionality reduction for moment-based distributionally robust optimization},
  journal       = {Operations Research},
  volume        = {74},
  number        = {3},
  pages         = {1648--1669},
  publisher     = {INFORMS},
  year          = {2026},
}

@inproceedings{paty2019subspace,
  author        = {Paty, Fran{\c{c}}ois-Pierre and Cuturi, Marco},
  title         = {Subspace robust {Wasserstein} distances},
  booktitle     = {Proceedings of the 36th International Conference on Machine Learning},
  series        = {Proceedings of Machine Learning Research},
  volume        = {97},
  pages         = {5072--5081},
  publisher     = {PMLR},
  year          = {2019},
}

@book{rockafellar1970,
  author        = {R. Tyrrell Rockafellar},
  title         = {Convex Analysis},
  series        = {Princeton Mathematical Series},
  volume        = {28},
  publisher     = {Princeton University Press},
  address       = {Princeton, NJ},
  year          = {1970},
}

@inproceedings{sinha2018,
  author        = {Aman Sinha and Hongseok Namkoong and John Duchi},
  title         = {Certifying some distributional robustness with principled adversarial training},
  booktitle     = {International Conference on Learning Representations},
  year          = {2018},
}

@article{xu2018,
  author        = {Huifu Xu and Yongchao Liu and Hailin Sun},
  title         = {Distributionally robust optimization with matrix moment constraints: {Lagrange} duality and cutting plane methods},
  journal       = {Mathematical Programming},
  volume        = {169},
  number        = {2},
  pages         = {489--529},
  year          = {2018},
}

@article{benTal2013,
  author        = {Ben-Tal, Aharon and {den Hertog}, Dick and {De Waegenaere}, Anja and Melenberg, Bertrand and Rennen, Gijs},
  title         = {Robust Solutions of Optimization Problems Affected by Uncertain Probabilities},
  journal       = {Management Science},
  volume        = {59},
  number        = {2},
  pages         = {341--357},
  year          = {2013},
}

@article{duchiNamkoong2021,
  author        = {Duchi, John C. and Namkoong, Hongseok},
  title         = {Learning Models with Uniform Performance via Distributionally Robust Optimization},
  journal       = {The Annals of Statistics},
  volume        = {49},
  number        = {3},
  pages         = {1378--1406},
  year          = {2021},
}

@article{xu2024flow,
  author        = {Xu, Chen and Lee, Jonghyeok and Cheng, Xiuyuan and Xie, Yao},
  title         = {Flow-based Distributionally Robust Optimization},
  journal       = {IEEE Journal on Selected Areas in Information Theory},
  volume        = {5},
  pages         = {62--77},
  year          = {2024},
}

@misc{zhu2024iterative,
  author        = {Zhu, Linglingzhi and Zhu, Yunqin and Xie, Yao},
  title         = {Distributionally robust optimization via iterative algorithms in continuous probability spaces},
  year          = {2024},
  howpublished  = {{arXiv} preprint arXiv:2412.20556},
  eprint        = {2412.20556},
  archivePrefix = {arXiv},
}

@article{cheng2018pca,
  author        = {Cheng, Jianqiang and Chen, Richard Li-Yang and Najm, Habib N. and Pinar, Ali and Safta, Cosmin and Watson, Jean-Paul},
  title         = {Distributionally Robust Optimization with Principal Component Analysis},
  journal       = {SIAM Journal on Optimization},
  volume        = {28},
  number        = {2},
  pages         = {1817--1841},
  year          = {2018},
}

@article{cheramin2022approximations,
  author        = {Cheramin, Meysam and Cheng, Jianqiang and Jiang, Ruiwei and Pan, Kai},
  title         = {Computationally Efficient Approximations for Distributionally Robust Optimization Under Moment and {Wasserstein} Ambiguity},
  journal       = {INFORMS Journal on Computing},
  volume        = {34},
  number        = {3},
  pages         = {1768--1794},
  year          = {2022},
}

@misc{wang2025representation,
  author        = {Wang, Zitao and Si, Nian and Liu, Molei},
  title         = {Representation-Aware Distributionally Robust Optimization: A Knowledge Transfer Framework},
  year          = {2025},
  howpublished  = {{arXiv} preprint arXiv:2509.09371},
  eprint        = {2509.09371},
  archivePrefix = {arXiv},
}

@misc{wu2026coreflow,
  author        = {Wu, Dongze and Zhu, Linglingzhi and Xie, Yao},
  title         = {{CoreFlow}: Low-Rank Matrix Generative Models},
  year          = {2026},
  howpublished  = {{arXiv} preprint arXiv:2604.24959},
  eprint        = {2604.24959},
  archivePrefix = {arXiv},
}

@misc{van2026distributionally,
  author        = {van Maarschalkerwaart, Floor and Mukherjee, Subhadip and Brune, Christoph and Carioni, Marcello},
  title         = {A Distributionally Robust Framework for Learned Reconstructions in Inverse Problems},
  year          = {2026},
  howpublished  = {{arXiv} preprint arXiv:2606.30230},
  eprint        = {2606.30230},
  archivePrefix = {arXiv},
}

@misc{xu2026distributionally,
  author        = {Xu, Chuang and Wood, Andrew T. A. and Yang, Yanrong},
  title         = {Distributionally Robust {PCA} with Data-Adaptive {Wasserstein} Geometry},
  year          = {2026},
  howpublished  = {{arXiv} preprint arXiv:2606.10463},
  eprint        = {2606.10463},
  archivePrefix = {arXiv},
}

@inproceedings{dohmatob2023origins,
  author        = {Dohmatob, Elvis and Guo, Chuan and Goibert, Morgane},
  title         = {Origins of Low-Dimensional Adversarial Perturbations},
  booktitle     = {Proceedings of the 26th International Conference on Artificial Intelligence and Statistics},
  series        = {Proceedings of Machine Learning Research},
  volume        = {206},
  pages         = {9221--9237},
  year          = {2023},
}

@misc{savostianova2024lowrank,
  author        = {Savostianova, Dayana and Zangrando, Emanuele and Tudisco, Francesco},
  title         = {Low-rank adversarial {PGD} attack},
  year          = {2024},
  howpublished  = {{arXiv} preprint arXiv:2410.12607},
  eprint        = {2410.12607},
  archivePrefix = {arXiv},
}

@misc{kazemi2021structured,
  author        = {Kazemi, Ehsan and Kerdreux, Thomas and Wang, Liquang},
  title         = {Generating Structured Adversarial Attacks Using {Frank--Wolfe} Method},
  year          = {2021},
  howpublished  = {{arXiv} preprint arXiv:2102.07360},
  eprint        = {2102.07360},
  archivePrefix = {arXiv},
}

@article{gui2026distributionally,
  author        = {Gui, Yu and Barber, Rina Foygel and Ma, Cong},
  title         = {Distributionally robust risk evaluation with an isotonic constraint},
  journal       = {Information and Inference: A Journal of the IMA},
  volume        = {15},
  number        = {3},
  pages         = {iaag025},
  publisher     = {Oxford University Press},
  year          = {2026},
}

@article{ye2005glram,
  author        = {Ye, Jieping},
  title         = {Generalized Low Rank Approximations of Matrices},
  journal       = {Machine Learning},
  volume        = {61},
  number        = {1--3},
  pages         = {167--191},
  year          = {2005},
}

@article{xie2014pmu,
  author        = {Xie, Le and Chen, Yang and Kumar, P. R.},
  title         = {Dimensionality reduction of synchrophasor data for early event detection: Linearized analysis},
  journal       = {IEEE Transactions on Power Systems},
  volume        = {29},
  number        = {6},
  pages         = {2784--2794},
  year          = {2014},
}

@article{gao2016pmu,
  author        = {Gao, Pengzhi and Wang, Meng and Ghiocel, Scott G. and Chow, Joe H. and Fardanesh, Bruce and Stefopoulos, George},
  title         = {Missing data recovery by exploiting low-dimensionality in power system synchrophasor measurements},
  journal       = {IEEE Transactions on Power Systems},
  volume        = {31},
  number        = {2},
  pages         = {1006--1013},
  year          = {2016},
}

@article{genton2007separable,
  author        = {Genton, Marc G.},
  title         = {Separable approximations of space-time covariance matrices},
  journal       = {Environmetrics},
  volume        = {18},
  number        = {7},
  pages         = {681--695},
  year          = {2007},
}

@article{yang2004twodpca,
  author        = {Yang, Jian and Zhang, David and Frangi, Alejandro F. and Yang, Jing-Yu},
  title         = {Two-dimensional {PCA}: A new approach to appearance-based face representation and recognition},
  journal       = {IEEE Transactions on Pattern Analysis and Machine Intelligence},
  volume        = {26},
  number        = {1},
  pages         = {131--137},
  year          = {2004},
}

@article{basri2003lambertian,
  author        = {Basri, Ronen and Jacobs, David W.},
  title         = {{Lambertian} Reflectance and Linear Subspaces},
  journal       = {IEEE Transactions on Pattern Analysis and Machine Intelligence},
  volume        = {25},
  number        = {2},
  pages         = {218--233},
  year          = {2003},
}

@article{mohajerin2018wasserstein,
  author        = {Mohajerin Esfahani, Peyman and Kuhn, Daniel},
  title         = {Data-Driven Distributionally Robust Optimization Using the {Wasserstein} Metric: Performance Guarantees and Tractable Reformulations},
  journal       = {Mathematical Programming},
  volume        = {171},
  number        = {1--2},
  pages         = {115--166},
  year          = {2018},
}

@article{mateos2013load,
  author        = {Mateos, Gonzalo and Giannakis, Georgios B.},
  title         = {Load curve data cleansing and imputation via sparsity and low rank},
  journal       = {IEEE Transactions on Smart Grid},
  volume        = {4},
  number        = {4},
  pages         = {2347--2355},
  publisher     = {IEEE},
  year          = {2013},
}

@article{wang2019factor,
  author        = {Wang, Dong and Liu, Xialu and Chen, Rong},
  title         = {Factor models for matrix-valued high-dimensional time series},
  journal       = {Journal of Econometrics},
  volume        = {208},
  number        = {1},
  pages         = {231--248},
  publisher     = {Elsevier},
  year          = {2019},
}

@article{blanchet2022optimal,
  title={Optimal transport-based distributionally robust optimization: Structural properties and iterative schemes},
  author={Blanchet, Jose and Murthy, Karthyek and Zhang, Fan},
  journal={Mathematics of Operations Research},
  volume={47},
  number={2},
  pages={1500--1529},
  year={2022},
  publisher={INFORMS}
}

@article{nietert2024robust,
  title={Robust distribution learning with local and global adversarial corruptions},
  author={Nietert, Sloan and Goldfeld, Ziv and Shafiee, Soroosh},
  journal={arXiv preprint arXiv:2406.06509},
  year={2024}
}

@book{rockafellar2009variational,
  title        = {Variational Analysis},
  author       = {Rockafellar, R. Tyrrell and Wets, Roger J.-B.},
  series       = {Grundlehren der Mathematischen Wissenschaften},
  volume       = {317},
  year         = {2009},
  publisher    = {Springer},
  address      = {Berlin, Heidelberg}
}
\end{document}